# Semi-harmonious and harmonious quasi-projection pairs on Hilbert $C^*$-modules

Xiaoyi Tian[a], Qingxiang Xu[a], Chunhong Fu[b]

[a]*Department of Mathematics, Shanghai Normal University, Shanghai 200234, PR China*
[b]*Health School Attached to Shanghai University of Medicine & Health Sciences, Shanghai 200237, PR China*

**Abstract**

For each adjointable idempotent $Q$ on a Hilbert $C^*$-module $H$, a specific projection $m(Q)$ called the matched projection of $Q$ was introduced recently due to the characterization of the minimum value among all the distances from projections to $Q$. Inspired by the relationship between $m(Q)$ and $Q$, another term called the quasi-projection pair $(P, Q)$ was also introduced recently, where $P$ is a projection on $H$ satisfying $Q^* = (2P - I)Q(2P - I)$, in which $Q^*$ is the adjoint operator of the idempotent $Q$ and $I$ is the identity operator on $H$. Some fundamental issues on quasi-projection pairs, such as the block matrix representations for quasi-projection pairs and the $C^*$-morphisms associated with quasi-projection pairs, are worthwhile to be investigated. This paper aims to make some detailed preparations. Two objects called the semi-harmonious quasi-projection pair and the harmonious quasi-projection pair are introduced and are systematically studied in the general setting of the adjointable operators on Hilbert $C^*$-modules. Some applications concerning the common similarity of operators and a norm equation associated with the Friedrichs angle are also dealt with. Furthermore, many examples are provided to illustrate the non-triviality of the associated characterizations.

*Email addresses:* tianxytian@163.com (Xiaoyi Tian), qingxiang_xu@126.com (Qingxiang Xu), fchlixue@163.com (Chunhong Fu)

## 1. Introduction and preliminaries

Let $H$ be a (right) Hilbert module over a $C^*$-algebra $\mathfrak{A}$, and $\mathcal{L}(H)$ be the set of all adjointable operators on $H$. For each $T \in \mathcal{L}(H)$, let $\mathcal{R}(T)$ and $\mathcal{N}(T)$ denote the range and the null space of $T$, respectively. An operator $Q \in \mathcal{L}(H)$ is called an idempotent if $Q^2 = Q$. If furthermore $Q$ is self-adjoint, then $Q$ is referred to be a projection. One fundamental issue concerning idempotents and projections is the characterization of the distances from all the projections to an idempotent $Q$, and some recent results on this issue can be found in [15, 16, 17, 21, 22]. For each idempotent $Q \in \mathcal{L}(H)$, a specific projection $m(Q)$ called the matched projection of $Q$ is introduced recently in [16]. It is proved in [22, Theorem 2.3, Corollary 2.5 and Remark 2.2] that

$$\|m(Q) - Q\| \le \|P - Q\| \le \|I - m(Q) - Q\|$$

for every projection $P \in \mathcal{L}(H)$, where $I$ denotes the identity operator on $H$. This shows that the distance away from $m(Q)$ to $Q$ (resp. from $I - m(Q)$ to $Q$) takes the minimum value (resp. the maximum value) among all the distances from projections to $Q$. The same is also true in the case that $H$ is a finite-dimensional Hilbert space and all the matrices are endued with the Frobenius norm [17, Theorems 3.2 and 4.2]. Inspired by the relationship between $m(Q)$ and $Q$, another term called the quasi-projection pair $(P, Q)$ is also introduced recently in [16], where $P \in \mathcal{L}(H)$ is a projection satisfying

$$Q^* = (2P - I)Q(2P - I), \tag{1.1}$$

in which $Q^*$ denotes the adjoint operator of the idempotent $Q$. In particular, $\big(m(Q), Q\big)$ is a quasi-projection pair, which is called the matched pair of $Q$.

Hilbert $C^*$-modules are the natural generalizations of Hilbert spaces by allowing the inner product to take values in certain $C^*$-algebra instead of the complex field [9, 13, 14]. It is known that a closed submodule of a Hilbert $C^*$-module may fail to be orthogonally complemented. Due to this dissimilarity between the Hilbert space and the Hilbert $C^*$-module, an adjointable operator on a Hilbert $C^*$-module may have no polar decomposition (see e.g. [10, Example 3.15]). So, some additional demanding is usually needed in dealing with a pair of projections on a general Hilbert $C^*$-module. Motivated by two norm equations used to characterize the Friedrichs angle [4], the term of the semi-harmonious pair $(P, Q)$ of projections is introduced in the sense that $\overline{\mathcal{R}(P + Q)}$ and $\overline{\mathcal{R}(2I - P - Q)}$ are both orthogonally complemented [6, Definition 1.1]. It is notable that such a semi-harmonious condition plays a crucial role in the study of the $C^*$-isomorphisms associated with two projections on a Hilbert $C^*$-module [6, Theorems 3.1 and 3.2]. Halmos' two projections theorem [2, 7] is well-known in the Hilbert space case as well as in the matrix case, which has many applications in various areas. Its generalization to the Hilbert $C^*$-module case is dealt with in [11, 20], and it is shown in [20, Theorem 3.3] that for every pair $(P,Q)$ of projections on a Hilbert $C^*$-module, the Halmos' two projections theorem is valid if and only if $(P,Q)$ is harmonious, namely, $(P,Q)$ and $(P, I-Q)$ are both semi-harmonious [11, Definition 4.1].

Although a quasi-projection pair $(P,Q)$ is composed of one projection and one idempotent rather than two projections, due to the restriction (1.1) it can be expected to obtain some results similar to the case of two projections mentioned above[1]. The key to achieving this goal is to clarify the exact meanings of the semi-harmony and the harmony for a quasi-projection pair separately, which is fulfilled in this paper; see Definitions 3.2 and 4.1 for the details.

A systematical study of semi-harmonious quasi-projection pairs is carried out in Section 3. To get a deeper understanding of semi-harmonious quasi-projection pairs, six operators and six closed submodules are defined at the beginning of this section for each quasi-projection pair. With the notations for these operators and submodules, three featured results are obtained in Section 3.1 for each quasi-projection pair. As a result, some equivalent conditions for a semi-harmonious quasi-projection pair are derived in Theorem 3.7, and are characterized furthermore in Corollaries 3.8, Lemma 3.10 and Theorem 3.12. A necessary and sufficient condition for a matched pair being semi-harmonious is provided in Theorem 3.9. A method to construct semi-harmonious quasi-projection pairs is revealed by Theorem 3.11. Two types of non-semi harmonious quasi-projection pairs are provided at the end of this section; see Examples 3.1 and 3.2 for the details.

Likewise, a systematical study of harmonious quasi-projection pairs is undertaken in Section 4. Some equivalent conditions for a harmonious quasi-projection pair are clarified in Theorem 4.2 and Corollary 4.3, and furthermore in Theorem 4.5 for some special quasi-projection pairs in a unital $C^*$-algebra, which serves as a Hilbert module over itself. The equivalence of a matched pair to be semi-harmonious and harmonious is revealed by Theorem 4.4. Attention has been paid to construct semi-harmonious quasi-projection pairs that are not harmonious (Examples 4.2 and 4.3).

As an application, the common similarity and the unitary equivalence of operators are dealt with in Section 5 in the setting of semi-harmonious quasi-projection pairs. For each semi-harmonious quasi-projection pair, an operator with non-zero parameters is defined by (5.12), and it is proved in Theorem 5.3 that such an operator is invertible and is in fact a common solution of (5.1) and (5.2). Employing (5.1) and (5.2) yields the modified equation (5.19), which leads obviously to equation (5.22). It is shown in Theorem 5.4 that for each semi-harmonious quasi-projection pair, equation (5.19) has a solution of unitary operator given by (5.21). Meanwhile, a counterexample is constructed in Theorem 5.5 to show that when the underlying quasi-projection pair fails to be semi-harmonious, equation (5.22) may have no solution of invertible operator. Another application is concerned with a norm equation associated with the Friedrichs angle, which is dealt with in Section 6 in the general setting of quasi-projection pairs. Under the weakest condition, the validity of equation (6.3) is verified in Theorem 6.4.

---

[1] Actually, some results on the block matrix representations for quasi-projection pairs have already been obtained in [15, Section 7].

We recall briefly some basic knowledge about adjointable operators on Hilbert $C^*$-modules. For more details, the reader is referred to [9, 13, 14]. Hilbert $C^*$-modules are generalizations of Hilbert spaces by allowing the inner product to take values in a $C^*$-algebra rather than in the complex field. Given Hilbert $C^*$-modules $H$ and $K$ over a $C^*$-algebra $\mathfrak{A}$, let $\mathcal{L}(H, K)$ denote the set of all adjointable operators from $H$ to $K$, with the abbreviation $\mathcal{L}(H)$ whenever $H = K$. Let $\mathcal{L}(H)_{\text{sa}}$ (resp. $\mathcal{L}(H)_+$) denote the set of all self-adjoint (resp. positive) elements in the $C^*$-algebra $\mathcal{L}(H)$. In the special case that $H$ is a Hilbert space, we use the notation $\mathbb{B}(H)$ instead of $\mathcal{L}(H)$.

Unless otherwise specified, throughout the rest of this paper $\mathbb{C}$ is the complex field, $\mathfrak{A}$ is a $C^*$-algebra, $E$, $H$ and $K$ are Hilbert $\mathfrak{A}$-modules. For every $C^*$-algebra $\mathfrak{B}$, $M_n(\mathfrak{B})$ denotes the $n \times n$ matrix algebra over $\mathfrak{B}$. Specifically, $M_n(\mathbb{C})$ stands for the set of all $n \times n$ complex matrices. The identity matrix and the zero matrix in $M_n(\mathbb{C})$ are denoted by $I_n$ and $0_n$, respectively.

It is known that every adjointable operator $A \in \mathcal{L}(H, K)$ is a bounded linear operator, which is also $\mathfrak{A}$-linear in the sense that

$$A(xa) = A(x)a, \quad \forall\, x \in H, a \in \mathfrak{A}. \tag{1.2}$$

Let $\mathcal{R}(A)$, $\mathcal{N}(A)$ and $A^*$ denote the range, the null space and the adjoint operator of $A$, respectively. When $H = K$, $A \in \mathcal{L}(H)$ and $M \subseteq H$, let $\sigma(A)$ denote the spectrum of $A$ with respect to $\mathcal{L}(H)$, and let $A|_M$ denote the restriction of $A$ on $M$. In the case that $A \in \mathcal{L}(H)_+$, the notation $A \geq 0$ is also used to indicate that $A$ is a positive operator on $H$.

Given a unital $C^*$-algebra $\mathfrak{B}$, its unit is denoted by $I_\mathfrak{B}$. When $V$ is a Hilbert module over a $C^*$-algebra $\mathfrak{C}$, the unit of $\mathcal{L}(V)$ will be denoted simply by $I_V$ rather than $I_{\mathcal{L}(V)}$. In most cases of this paper, the Hilbert $C^*$-module under consideration is represented by $H$, so the symbol $I_H$ is furthermore simplified to be $I$. In other words, the Italic letter $I$ is reserved for the identity operator on the Hilbert $C^*$-module $H$.

A closed submodule $M$ of $H$ is referred to be orthogonally complemented in $H$ if $H = M + M^\perp$, where
$$M^\perp = \{x \in H : \langle x, y \rangle = 0 \text{ for every } y \in M\}.$$

In this case, the projection from $H$ onto $M$ is denoted by $P_M$. We use the notations "$\oplus$" and "$\dotplus$" with different meanings. For Hilbert $\mathfrak{A}$-modules $H_1$ and $H_2$, let
$$H_1 \oplus H_2 = \left\{(h_1, h_2)^T : h_i \in H_i, i = 1, 2\right\},$$
which is also a Hilbert $\mathfrak{A}$-module whose $\mathfrak{A}$-valued inner product is given by
$$\langle (x_1, y_1)^T, (x_2, y_2)^T \rangle = \langle x_1, x_2 \rangle + \langle y_1, y_2 \rangle$$
for $x_i \in H_1, y_i \in H_2, i = 1, 2$. If both $H_1$ and $H_2$ are submodules of a Hilbert $\mathfrak{A}$-module $H$ such that $H_1 \cap H_2 = \{0\}$, then we set
$$H_1 \dotplus H_2 = \{h_1 + h_2 : h_i \in H_i, i = 1, 2\}.$$

For each $T \in \mathcal{L}(H, K)$, the notation $|T|$ is used to denote the square root of $T^*T$. Hence, $|T^*| = (TT^*)^{\frac{1}{2}}$. From [10, Lemma 3.9] we know that there exists at most a partial isometry $U \in \mathcal{L}(H, K)$ satisfying
$$T = U|T| \quad \text{with} \quad U^*U = P_{\overline{\mathcal{R}(T^*)}}. \tag{1.3}$$

If such a partial isometry $U$ exists, then the representation (1.3) is called the polar decomposition of $T$ [10, Definition 3.10]. As shown in [10, Example 3.15], an adjointable operator on a Hilbert $C^*$-module may fail to have the polar decomposition. The following lemma gives a characterization of the existence of the polar decomposition.

**Lemma 1.1.** [10, Lemma 3.6 and Theorem 3.8] *For every $T \in \mathcal{L}(H, K)$, the following statements are equivalent:*

(i) *$T$ has the polar decomposition;*
(ii) *$T^*$ has the polar decomposition;*
(iii) *$\overline{\mathcal{R}(T)}$ and $\overline{\mathcal{R}(T^*)}$ are orthogonally complemented in $K$ and $H$, respectively.*

*If* (i)–(iii) *are satisfied, and the polar decomposition of $T$ is given by* (1.3), *then the polar decomposition of $T^*$ is represented by*

$$T^* = U^*|T^*| \quad \text{with} \quad UU^* = P_{\overline{\mathcal{R}(T)}}.$$

*Remark* 1.1. Suppose that $T \in \mathcal{L}(H,K)$ has the polar decomposition represented by (1.3). For simplicity, henceforth we just say that $T$ has the polar decomposition $T = U|T|$.

To achieve the main results of this paper, we need some additional known lemmas.

**Lemma 1.2.** [9, Proposition 3.7] *For every $T \in \mathcal{L}(H,K)$, $\overline{\mathcal{R}(T^*T)} = \overline{\mathcal{R}(T^*)}$ and $\overline{\mathcal{R}(TT^*)} = \overline{\mathcal{R}(T)}$.*

**Lemma 1.3.** ([10, Proposition 2.9] and [18, Lemma 2.2]) *Let $T \in \mathcal{L}(H)_+$. Then for every $\alpha > 0$, $\overline{\mathcal{R}(T^\alpha)} = \overline{\mathcal{R}(T)}$ and $\mathcal{N}(T^\alpha) = \mathcal{N}(T)$.*

**Lemma 1.4.** [10, Proposition 2.7] *Let $A \in \mathcal{L}(H,K)$ and $B, C \in \mathcal{L}(E,H)$ be such that $\overline{\mathcal{R}(B)} = \overline{\mathcal{R}(C)}$. Then $\overline{\mathcal{R}(AB)} = \overline{\mathcal{R}(AC)}$.*

The rest of the paper is organized as follows. In Section 2, some basic knowledge about the quasi-projection pair and the matched projection are provided, and three examples of quasi-projection pairs are newly added. Sections 3 and 4 are devoted to the study of semi-harmonious quasi-projection pairs and harmonious quasi-projection pairs, respectively. Some applications concerning the common similarity of operators and a norm equation associated with the Friedrichs angle are dealt with in Sections 5 and 6, respectively.

## 2. Characterizations and examples of quasi-projection pairs

In this section, we are concerned with the quasi-projection pair, which is introduced recently in [16] as follows.

**Definition 2.1.** [16, Definition 1.2] An ordered pair $(P,Q)$ is called a quasi-projection pair on $H$ if $P \in \mathcal{L}(H)$ is a projection, while $Q \in \mathcal{L}(H)$ is an idempotent such that

$$PQ^*P = PQP, \quad PQ^*(I-P) = -PQ(I-P), \quad (2.1)$$
$$(I-P)Q^*(I-P) = (I-P)Q(I-P). \quad (2.2)$$

*Remark* 2.1. [16, Remark 1.3] From Definition 2.1, it is clear that $(P, P)$ is a quasi-projection pair whenever $P \in \mathcal{L}(H)$ is a projection. In the special case that $P$ is a trivial projection (that is, either $P = I$ or $P = 0$), then for every $Q \in \mathcal{L}(H)$, $(P, Q)$ is a quasi-projection pair if and only if $Q$ is a projection.

Two useful characterizations of quasi-projection pairs can be stated as follows.

**Lemma 2.1.** [16, Theorem 1.4] Suppose that $P \in \mathcal{L}(H)$ is a projection and $Q \in \mathcal{L}(H)$ is an idempotent. If one element in

$$\Sigma := \{(A, B) : A \in \{P, I - P\}, B \in \{Q, Q^*, I - Q, I - Q^*\}\}$$

is a quasi-projection pair, then all the remaining elements in $\Sigma$ are quasi-projection pairs.

**Lemma 2.2.** [16, Theorem 1.5] Suppose that $P \in \mathcal{L}(H)$ is a projection and $Q \in \mathcal{L}(H)$ is an idempotent. Then the following statements are equivalent:

(i) $(P, Q)$ is a quasi-projection pair;

(ii) $Q^* = (2P - I)Q(2P - I)$;

(iii) $|Q^*| = (2P - I)|Q|(2P - I)$.

Now, we turn to the description of the matched projection, and recall two properties about the matched projection.

**Definition 2.2.** [16, Definition 2.2 and Theorem 2.4] For each idempotent $Q \in \mathcal{L}(H)$, its matched projection, written $m(Q)$, is defined by

$$m(Q) = \frac{1}{2}(|Q^*| + Q^*)|Q^*|^\dagger (|Q^*| + I_n)^{-1}(|Q^*| + Q),$$

where $|Q^*|^\dagger$ denotes the Moore-Penrose inverse of $|Q^*|$, which can be expressed as

$$|Q^*|^\dagger = \left(P_{\mathcal{R}(Q)} P_{\mathcal{R}(Q^*)} P_{\mathcal{R}(Q)}\right)^{\frac{1}{2}}.$$

In the rest of this paper, the notation $m(Q)$ is reserved to denote the matched projection of an idempotent $Q$.

**Definition 2.3.** For each idempotent $Q \in \mathcal{L}(H)$, $(m(Q), Q)$ is called the matched pair of $Q$ on $H$.

*Remark* 2.2. It is remarkable that for every idempotent $Q \in \mathcal{L}(H)$, its matched pair is a quasi-projection pair [16, Theorem 2.1].

**Lemma 2.3.** [16, Remark 2.10] *For every idempotent $Q \in \mathcal{L}(H)$, we have*
$$\mathcal{R}\big[m(Q)\big] = \mathcal{R}(|Q^*| + Q^*) = \mathcal{R}(|Q| + Q).$$

*Remark* 2.3. It is notable that $\mathcal{R}(|Q^*|) = \mathcal{R}(Q)$ for every idempotent $Q$ on $H$. Indeed, it is obvious if $Q$ is a projection. When $Q$ is not a projection, this equality can be derived by using the standard $2 \times 2$ block matrix representation of $Q$. In what follows, such an equality will be used without specified.

**Lemma 2.4.** [16, Theorems 2.7 and 2.14] *For every idempotent $Q \in \mathcal{L}(H)$, we have*
$$m(Q^*) = m(Q), \quad I - m(Q) = m(I - Q).$$

It is worthwhile to find out some examples of quasi-projection pairs. For this, three examples are newly provided as follows.

**Example 2.1.** Recall that an operator $J \in \mathcal{L}(H)$ is referred to be a symmetry if $J^2 = I$ and $J = J^*$. As in the Hilbert space case [1], we call such a pair $(H, J)$ (or simply $H$) as a Krein space. An operator $Q \in \mathcal{L}(H)$ is said to be a weighted projection (on the indefinite inner-product space induced by $J$) if $(JQ)^* = JQ$.

Let $J_+$ be the positive part of $J$ given by $J_+ = \frac{1}{2}(I + J)$, and $Q \in \mathcal{L}(H)$ be an arbitrary idempotent. It is clear that $J_+$ is a projection, and by Lemma 2.2(ii) $(J_+, Q)$ is a quasi-projection pair if and only if
$$Q^* = (2J_+ - I)Q(2J_+ - I),$$
or equivalently $Q^* = JQJ$, which is exactly the case that $(JQ)^* = JQ$, since $J$ is a symmetry. This shows that $(J_+, Q)$ is a quasi-projection pair if and only if $Q$ is a weighted projection. Let $J_- = I - J_+$, which is called the negative part of $J$. From Lemma 2.1, we see that for every idempotent $Q \in \mathcal{L}(H)$, $(J_+, Q)$ is a quasi-projection if and only $(J_-, Q)$ is.

Given a projection $P$, there might exist many idempotents $Q$ that make $(P, Q)$ to be a quasi-projection pair. This will be shown in our next example.

**Example 2.2.** Suppose that $A \in \mathcal{L}(H)_{\text{sa}}$ satisfying $A^2 - A \geq 0$. Let
$$\ell(A) = (A^2 - A)^{\frac{1}{2}}, \tag{2.3}$$

and let $P, Q \in \mathcal{L}(H \oplus H)$ be given by

$$P = \begin{pmatrix} I & 0 \\ 0 & 0 \end{pmatrix}, \quad Q = \begin{pmatrix} A & -\ell(A) \\ \ell(A) & I - A \end{pmatrix}.$$

It is a routine matter to show that $Q$ is an idempotent, and both of (2.1) and (2.2) are satisfied. So, $(P, Q)$ defined as above is a quasi-projection pair.

The following example shows that given an idempotent $Q$, there might have only finite projections $P$ to make $(P, Q)$ as a quasi-projection pair.

**Example 2.3.** Let $Q = \begin{pmatrix} 1 & a \\ 0 & 0 \end{pmatrix}$ with $a \in \mathbb{C} \setminus \{0\}$. It can be shown (see e.g. [15, Example 3.10]) that there exist only four projections $P_i \in M_2(\mathbb{C})$ such that $(P_i, Q)$ $(0 \leq i \leq 3)$ are quasi-projection pairs, where $b = \sqrt{1 + |a|^2}$ and

$$P_0 = m(Q) = \frac{1}{2b} \begin{pmatrix} b+1 & a \\ \bar{a} & b-1 \end{pmatrix}, \quad P_1 = \frac{1}{2b} \begin{pmatrix} b-1 & a \\ \bar{a} & b+1 \end{pmatrix},$$

while $P_2 = I_2 - P_0$ and $P_3 = I_2 - P_1$.

## 3. Semi-harmonious quasi-projection pairs

In this section, we will introduce and study semi-harmonious quasi-projection pairs. It is helpful to specify some frequently used symbols for those operators and submodules associated with two idempotents. On the whole, we will introduce six operators and six submodules as follows.

**Definition 3.1.** Given idempotents $P, Q \in \mathcal{L}(H)$, let

$$T_1 = P(I - Q), \quad T_2 = (I - P)Q, \tag{3.1}$$
$$T_3 = PQ(I - P), \quad T_4 = (I - P)QP, \tag{3.2}$$
$$\widetilde{T_1} = T_1(2P - 1), \quad \widetilde{T_2} = -T_2(2P - I), \tag{3.3}$$

and let

$$H_1 = \mathcal{R}(P) \cap \mathcal{R}(Q), \quad H_4 = \mathcal{N}(P) \cap \mathcal{N}(Q), \tag{3.4}$$
$$H_2 = \mathcal{R}(P) \cap \mathcal{N}(Q), \quad H_3 = \mathcal{N}(P) \cap \mathcal{R}(Q), \tag{3.5}$$
$$H_5 = \overline{\mathcal{R}(T_3)}, \quad H_6 = \overline{\mathcal{R}(T_4)}. \tag{3.6}$$

### 3.1. Featured results for every quasi-projection pair

We begin with a characterization of the adjoint operators as follows.

**Lemma 3.1.** *Given a quasi-projection pair $(P,Q)$ on $H$, let $T_i$ ($1 \le i \le 4$) and $\widetilde{T}_i$ ($1 \le i \le 2$) be defined by (3.1)–(3.3). Then*

$$T_1^* = (2P - I)(I - Q)P, \quad \widetilde{T}_1^{\,*} = (I - Q)P, \qquad (3.7)$$
$$T_2^* = -(2P - I)Q(I - P), \quad \widetilde{T}_2^{\,*} = Q(I - P), \qquad (3.8)$$
$$T_3^* = -T_4, \quad T_4^* = -T_3. \qquad (3.9)$$

*Proof.* Since $P$ ia a projection, we know that $2P - I$ is a self-adjoint unitary (symmetry). Hence, by Lemma 2.2(ii) we have

$$T_2^* = (2P - I)Q(2P - I) \cdot (I - P) = -(2P - I)Q(I - P),$$
$$\widetilde{T}_2^{\,*} = -(2P - I)T_2^* = Q(I - P),$$
$$T_4^* = P \cdot (2P - I)Q(2P - I) \cdot (I - P) = -T_3.$$

This gives the verification of (3.8) and the second equality in (3.9). In view of Lemma 2.1, the rest equalities can be obtained by a simple replacement of $(P,Q)$ with $(I - P, I - Q)$. □

Next, we provide three featured results for every quasi-projection pair.

**Lemma 3.2.** *For every quasi-projection pair $(P,Q)$ on $H$, $T_1 + T_2$ is self-adjoint, where $T_1$ and $T_2$ are defined by (3.1).*

*Proof.* From (3.7) and (3.8), we have

$$(T_1 + T_2)^* = (2P - I)\big[(I - Q)P - Q(I - P)\big] = P - 2PQ + Q$$
$$= P(I - Q) + (I - P)Q = T_1 + T_2.$$
□

**Lemma 3.3.** *For every quasi-projection pair $(P,Q)$ on $H$, we have*

$$|T_i| \cdot |T_j^*|^2 = |T_i|^2 \cdot |T_j^*|, \quad i,j = 1,2,$$

*where $T_1$ and $T_2$ are defined by (3.1).*

*Proof.* Given an arbitrary operator $T \in \mathcal{L}(H)$, it is clear that $|T^*|^2 \cdot T = T \cdot |T|^2$, which leads to

$$p(|T^*|^2) \cdot T = T \cdot p(|T|^2)$$

for every polynomial $p$. Hence,
$$|T^*| \cdot T = T \cdot |T| \tag{3.10}$$
by the functional calculus of elements in a commutative $C^*$-algebra, since the function $t \to \sqrt{t}$ can be approximated uniformly by a sequence of polynomials on any compact subset of $[0, +\infty)$. Also, from Lemmas 1.3 and 1.2 we can obtain
$$\overline{\mathcal{R}(|T^*|)} = \overline{\mathcal{R}(T)}, \quad \mathcal{N}(|T|) = \mathcal{N}(T). \tag{3.11}$$
In virtue of (3.11) and (3.10), we see that for any $T, S \in \mathcal{L}(H)$,
$$\begin{aligned} |T| \cdot |S^*|^2 = |T|^2 \cdot |S^*| &\Longleftrightarrow |T|\big(|S^*| - |T|\big)|S^*| = 0 \\ &\Longleftrightarrow |T|\big(|S^*| - |T|\big)S = 0 \\ &\Longleftrightarrow T\big(|S^*| - |T|\big)S = 0 \\ &\Longleftrightarrow T \cdot |S^*| \cdot S = T \cdot |T| \cdot S \\ &\Longleftrightarrow T \cdot |S^*| \cdot S = |T^*| \cdot T \cdot S. \end{aligned}$$

So, our aim here is to prove that
$$T_i \cdot |T_j^*| \cdot T_j = |T_i^*| \cdot T_i \cdot T_j \quad (i, j = 1, 2). \tag{3.12}$$

First, we consider the case that $i = j = 1$. From (3.1), we have
$$T_1 P = P(I - Q)P,$$
which is self-adjoint according to (2.1). Utilizing the expression of $T_1^*$ given in (3.7) yields
$$T_1 T_1^* = P(I - Q) \cdot (2P - I)(I - Q)P = 2(T_1 P)^2 - T_1 P.$$
Hence,
$$T_1 P \cdot T_1 T_1^* = 2(T_1 P)^3 - (T_1 P)^2.$$
As $T_1 P = (T_1 P)^*$, this shows that
$$T_1 P \cdot T_1 T_1^* = (T_1 P \cdot T_1 T_1^*)^* = T_1 T_1^* \cdot T_1 P.$$
Similar reasoning in the derivation of (3.10) gives
$$T_1 P \cdot |T_1^*| = |T_1^*| \cdot T_1 P,$$
which can be simplified as
$$T_1 \cdot |T_1^*| = |T_1^*| \cdot T_1 P, \tag{3.13}$$

since $P$ is a projection, and by (3.11) and (3.1) $\overline{\mathcal{R}(|T_1^*|)} = \overline{\mathcal{R}(T_1)} \subseteq \mathcal{R}(P)$. It follows that

$$T_1 \cdot |T_1^*| \cdot T_1 = |T_1^*| \cdot T_1 P \cdot T_1 = |T_1^*| \cdot T_1^2$$

by post-multiplying $T_1$ on both sides of (3.13). So, the validity of (3.12) is verified in the case considered as above.

Next, we show that (3.12) is true in the case that $i = 1$ and $j = 2$. A simple calculation gives

$$T_1(I - P) = -PQ(I - P), \quad T_1 PQ = -T_1 T_2. \tag{3.14}$$

From (3.7) and the observation $T_1 Q = 0$, we have

$$T_1 T_1^* = T_1 \cdot (2P - I)(I - Q)P = T_1 \cdot (P - 2PQP).$$

Combining the above expression for $T_1 T_1^*$ with (3.14) and Lemma 2.2(ii) yields

$$\begin{aligned} T_1 T_1^* \cdot T_1(I - P) &= -T_1(P - 2PQP) \cdot PQ(I - P) \\ &= -T_1 PQ \cdot (I - 2P)Q(I - P) \\ &= T_1 T_2 \cdot Q^*(I - P). \end{aligned}$$

Hence,

$$T_1 T_1^* \cdot T_1(I - P) = T_1 T_2 T_2^*,$$

which means that

$$T_1 T_1^* \cdot T_1(I - P) = T_1(I - P) \cdot T_2 T_2^*,$$

since $(I - P)T_2 = T_2$. Consequently,

$$|T_1^*| \cdot T_1(I - P) = T_1(I - P) \cdot |T_2^*|,$$

which leads by the observation $\overline{\mathcal{R}(|T_2^*|)} = \overline{\mathcal{R}(T_2)} \subseteq \mathcal{R}(I - P)$ to

$$|T_1^*| \cdot T_1(I - P) = T_1 \cdot |T_2^*|.$$

Therefore,

$$T_1 \cdot |T_2^*| \cdot T_2 = |T_1^*| \cdot T_1(I - P) \cdot T_2 = |T_1^*| \cdot T_1 \cdot T_2.$$

So, the desired conclusion follows.

Finally, we may replace $(P, Q)$ with $(I - P, I - Q)$ to get (3.12) in the remaining two cases. □

**Lemma 3.4.** *For every quasi-projection pair $(P,Q)$ on $H$, we have*

$$\begin{aligned} H_1 &= \mathcal{R}(P) \cap \mathcal{R}(Q^*), & H_4 &= \mathcal{N}(P) \cap \mathcal{N}(Q^*), \\ H_2 &= \mathcal{R}(P) \cap \mathcal{N}(Q^*), & H_3 &= \mathcal{N}(P) \cap \mathcal{R}(Q^*), \end{aligned} \quad (3.15)$$

*where $H_i$ $(1 \le i \le 4)$ are defined by (3.4)–(3.5).*

*Proof.* Let $x \in \mathcal{R}(P)$. It is obvious that $(2P - I)x = x$. If furthermore $x \in \mathcal{R}(Q)$, then $Qx = x$ and from Lemma 2.2(ii) we have

$$Q^*x = (2P - I)Q(2P - I)x = x.$$

It follows that $H_1 \subseteq \mathcal{R}(P) \cap \mathcal{R}(Q^*)$. Since $2P - I$ is a symmetry, similar reasoning gives $\mathcal{R}(P) \cap \mathcal{R}(Q^*) \subseteq H_1$. This completes the verification of the first equation in (3.15).

Replacing $(P,Q)$ with $(I-P, I-Q)$, $(P, I-Q)$ and $(I-P, Q)$ respectively, the remaining three equations can be derived immediately. □

*3.2. Equivalent conditions for the semi-harmonious quasi-projection pair*

As shown in [10, Example 3.15], an adjointable operator on a Hilbert $C^*$-module may have no polar decomposition. This leads us to make a definition as follows.

**Definition 3.2.** Let $(P,Q)$ be a quasi-projection pair. If $T_1$ and $T_2$ defined by (3.1) both have the polar decompositions, then $(P,Q)$ is said to be semi-harmonious.

Suppose that $(P,Q)$ is a quasi-projection pair on $H$. In Theorem 3.7, we will prove that the above demanding for $(P,Q)$ being semi-harmonious can in fact be weakened as one of $T_1$ and $T_2$ has the polar decomposition. To get such a characterization, we need two useful lemmas as follows.

**Lemma 3.5.** *For every idempotents $P, Q \in \mathcal{L}(H)$, the following statements are valid:*

(i) $\mathcal{N}[(I-Q)(I-P)] = \mathcal{R}(P) + H_3$;
(ii) $\mathcal{N}(T_4) = \mathcal{N}(P) + H_1 + H_2$,

*where $H_1, H_2, H_3$ and $T_4$ are defined by (3.4), (3.5) and (3.2), respectively.*

*Proof.* A proof given in [5, Lemma 2.1] for two projections is also valid in the case that $P, Q \in \mathcal{L}(H)$ are idempotents. □

**Lemma 3.6.** *For every quasi-projection pair $(P,Q)$ on $H$, let $T_1, T_2, \widetilde{T_1}, \widetilde{T_2}$ and $H_1$ be defined by (3.1), (3.3) and (3.4), respectively. Then the following statements are equivalent:*

  (i) $\overline{\mathcal{R}(\widetilde{T_1})}$ *is orthogonally complemented in $H$;*
  (ii) $\overline{\mathcal{R}(T_1)}$ *is orthogonally complemented in $H$;*
  (iii) $\overline{\mathcal{R}(\widetilde{T_2}^*)}$ *is orthogonally complemented in $H$;*
  (iv) $\overline{\mathcal{R}(T_2^*)}$ *is orthogonally complemented in $H$;*
  (v) $\overline{\mathcal{R}\big[I - P + (I-Q)(I-Q)^*\big]}$ *is orthogonally complemented in $H$.*

*In each case, $H_1$ is orthogonally complemented in $H$ such that*

$$P = P_{\overline{\mathcal{R}(T_1)}} + P_{H_1}, \quad P_{\mathcal{R}(Q)} = P_{\overline{\mathcal{R}(\widetilde{T_2}^*)}} + P_{H_1}. \tag{3.16}$$

*Proof.* By the first equation in (3.15) and the second equation in (3.8), we have

$$\overline{\mathcal{R}(T_1)} \perp H_1 \quad \text{and} \quad \overline{\mathcal{R}(\widetilde{T_2}^*)} \perp H_1.$$

(i)⟺(ii) and (iii)⟺(iv). From (3.3) and the invertibility of $2P - I$, it is obvious that $\mathcal{R}(\widetilde{T_1}) = \mathcal{R}(T_1)$, and thus (i)⟺(ii) follows. From (3.8) we have $T_2^* = (I - 2P)\widetilde{T_2}^*$, which clearly leads to (iii)⟺(iv), since $I - 2P$ is a unitary.

(v)⟹(ii) and (v)⟹(iii). By assumption, we have $H = \overline{\mathcal{R}(S)} \dotplus \mathcal{N}(S)$, where $S$ is the positive operator given by

$$S = I - P + (I - Q)(I - Q)^*. \tag{3.17}$$

Clearly, $\mathcal{N}(S) = \mathcal{R}(P) \cap \mathcal{R}(Q^*)$, which is equal to $H_1$ by (3.15). Hence, $H_1$ is orthogonally complemented in $H$ such that

$$H = \overline{\mathcal{R}(S)} \dotplus H_1. \tag{3.18}$$

So each $x \in H$ can be decomposed as $x = x_1 + x_2$ for some $x_1 \in \overline{\mathcal{R}(S)}$ and $x_2 \in H_1$, which gives

$$Px = Px_1 + x_2. \tag{3.19}$$

Note that

$$P\overline{\mathcal{R}(S)} \subseteq \overline{\mathcal{R}(PS)} = \overline{\mathcal{R}\big[(T_1(I-Q)^*\big]} \subseteq \overline{\mathcal{R}(T_1)},$$

so by (3.19) we have

$$\mathcal{R}(P) \subseteq \overline{\mathcal{R}(T_1)} + H_1 \subseteq \mathcal{R}(P),$$

which leads to the orthogonal decomposition

$$\mathcal{R}(P) = \overline{\mathcal{R}(T_1)} \dotplus H_1. \tag{3.20}$$

Consequently, we arrive at

$$H = \overline{\mathcal{R}(T_1)} \dotplus \overline{\mathcal{R}(T_1)}^\perp,$$

where

$$\overline{\mathcal{R}(T_1)}^\perp = \mathcal{R}(P)^\perp + H_1.$$

Utilizing (3.20) yields the first equation in (3.16). Similarly, from (3.18) and the second equation in (3.8), it can be obtained that

$$\mathcal{R}(Q) = \overline{\mathcal{R}(\widetilde{T_2}^*)} + H_1.$$

So, $\overline{\mathcal{R}(\widetilde{T_2}^*)}$ is orthogonally complemented in $H$ such that the second equation in (3.16) is satisfied.

(iv)$\Longrightarrow$(ii). By assumption, we have

$$PH = P\overline{\mathcal{R}(T_2^*)} + P\mathcal{N}(T_2). \tag{3.21}$$

A direct application of Lemma 3.5(i) gives $\mathcal{N}(T_2) = \mathcal{R}(I - Q) + H_1$, which yields

$$P\mathcal{N}(T_2) = \mathcal{R}(T_1) + H_1. \tag{3.22}$$

In addition, from the first equation in (3.8) it is easily seen that

$$PT_2^* = T_1(I - P).$$

Hence, the orthogonal decomposition (3.20) can be derived by combining the above equation with (3.21) and (3.22). Therefore, $\overline{\mathcal{R}(T_1)}$ is orthogonally complemented in $H$.

(ii)$\Longrightarrow$(v). By assumption, we have $H = \overline{\mathcal{R}(T_1)} + \mathcal{N}(T_1^*)$. From the first equation in (3.7) together with the invertibility of $2P - I$ and Lemma 3.5(i), it can be deduced that

$$\mathcal{N}(T_1^*) = \mathcal{N}\big[(I - Q)P\big] = \mathcal{R}(I - P) + H_1.$$

This shows that

$$H = \overline{\mathcal{R}(T_1)} + \mathcal{R}(I - P) + H_1. \tag{3.23}$$

Meanwhile, from Lemma 1.2 and [5, Lemma 3.5] we can obtain

$$\overline{\mathcal{R}(I-P)+\mathcal{R}(I-Q)} = \overline{\mathcal{R}(I-P)+\mathcal{R}\big[(I-Q)(I-Q)^*\big]} = \overline{\mathcal{R}(S)},$$

where $S$ is defined by (3.17). It follows that

$$\mathcal{R}(I-P) \subseteq \overline{\mathcal{R}(S)} \quad \text{and} \quad \overline{\mathcal{R}(I-Q)} \subseteq \overline{\mathcal{R}(S)},$$

which mean that

$$\mathcal{R}(T_1) \subseteq \mathcal{R}(I-P) + \mathcal{R}(I-Q) \subseteq \overline{\mathcal{R}(S)}, \tag{3.24}$$

since $T_1$ can be expressed alternatively as

$$T_1 = -(I-P)(I-Q) + I - Q.$$

Combining (3.23) and (3.24) yields the orthogonal decomposition (3.18). □

Now, we are in the position to provide some equivalent conditions for the semi-harmonious quasi-projection pair.

**Theorem 3.7.** *For every quasi-projection pair $(P,Q)$ on $H$, let $T_1$, $T_2$, $\widetilde{T_1}$, $\widetilde{T_2}$, $H_1$ and $H_4$ be defined by (3.1), (3.3) and (3.4), respectively. Then the following statements are equivalent:*

 (i) *$(P,Q)$ is semi-harmonious;*
 (ii) *$\overline{\mathcal{R}(T_i)}(i=1,2)$ are both orthogonally complemented in $H$;*
 (iii) *$\overline{\mathcal{R}(T_i^*)}(i=1,2)$ are both orthogonally complemented in $H$;*
 (iv) *$\overline{\mathcal{R}(\widetilde{T_i})}(i=1,2)$ are both orthogonally complemented in $H$;*
 (v) *$\overline{\mathcal{R}(\widetilde{T_i}^*)}(i=1,2)$ are both orthogonally complemented in $H$;*
 (vi) *$\overline{\mathcal{R}(T_1)}$ and $\overline{\mathcal{R}(T_1^*)}$ are both orthogonally complemented in $H$;*
 (vii) *$T_1$ has the polar decomposition;*
 (viii) *$\overline{\mathcal{R}(T_1^*)}$ and $\overline{\mathcal{R}(\widetilde{T_2}^*)}$ are both orthogonally complemented in $H$;*
 (ix) *$\overline{\mathcal{R}(T_2)}$ and $\overline{\mathcal{R}(T_2^*)}$ are both orthogonally complemented in $H$;*
 (x) *$T_2$ has the polar decomposition;*
 (xi) *$\overline{\mathcal{R}(T_2^*)}$ and $\overline{\mathcal{R}(\widetilde{T_1}^*)}$ are both orthogonally complemented in $H$.*

*In each case, $H_1$ and $H_4$ are orthogonally complemented in $H$ such that (3.16) is satisfied, and*

$$I - P = P_{\overline{\mathcal{R}(T_2)}} + P_{H_4}, \quad P_{\mathcal{N}(Q)} = P_{\overline{\mathcal{R}(\widetilde{T_1}^*)}} + P_{H_4}. \tag{3.25}$$

*Proof.* The conclusion follows directly from Lemmas 1.1 and 3.6, together with the replacement of $(P,Q)$ by $(I-P, I-Q)$. □

**Corollary 3.8.** *For every quasi-projection pair $(P,Q)$ on $H$, if one pair in*

$$(P,Q), \quad (I-P, I-Q), \quad (P, Q^*) \quad \text{and} \quad (I-P, I-Q^*) \tag{3.26}$$

*is semi-harmonious, then the remaining three pairs are all semi-harmonious.*

*Proof.* In view of Lemma 2.1, each pair in (3.26) is a quasi-projection pair. From the definition of $T_1$ and $T_2$ together with (i)⟺(ii) in Theorem 3.7, it follows easily that

$(P,Q)$ is semi-harmonious $\iff$ $(I-P, I-Q)$ is semi-harmonious,
$(P, Q^*)$ is semi-harmonious $\iff$ $(I-P, I-Q^*)$ is semi-harmonious.

Here the point is,

$$\mathcal{R}(T_1) = \mathcal{R}\big[P(I-Q)^*\big], \quad \mathcal{R}(T_2) = \mathcal{R}\big[(I-P)Q^*\big]. \tag{3.27}$$

Indeed, by Lemma 2.2(ii) we have

$$(I-P)Q^* = (I-P)(2P-I)Q(2P-I) = T_2(I-2P),$$

which leads by the invertibility of $I-2P$ to the second equation in (3.27). Replacing $(P,Q)$ with $(I-P, I-Q)$ yields the first equation in (3.27). Hence, $(P,Q)$ is semi-harmonious if and only if $(P, Q^*)$ is semi-harmonious. This completes the proof. □

For the matched pair, we have a specific characterization as follows.

**Theorem 3.9.** *For every idempotent $Q \in \mathcal{L}(H)$, its matched pair is semi-harmonious if and only if $P_{\mathcal{R}(Q)} - Q$ has the polar decomposition.*

*Proof.* What we are concerned with is the quasi-projection pair $(P,Q)$ with $P = m(Q)$. So, by (3.7), (3.8) and Lemma 2.4 we have

$$T_1^* = (I - Q^*) m(Q), \quad \widetilde{T_2}^* = Q\big[I - m(Q)\big] = Q \cdot m(I-Q). \tag{3.28}$$

Since $\mathcal{R}\big[m(Q)\big] = \mathcal{R}(|Q^*| + Q^*)$ (see Lemma 2.3) and $\mathcal{R}(|Q^*|) = \mathcal{R}\big[P_{\mathcal{R}(Q)}\big]$, it can be deduced by Lemma 1.4 that

$$\overline{\mathcal{R}(T_1^*)} = \overline{\mathcal{R}\big[(I-Q^*)(|Q^*| + Q^*)\big]} = \overline{\mathcal{R}\big[(I-Q^*)|Q^*|\big]}$$
$$= \overline{\mathcal{R}\big[(I-Q^*) P_{\mathcal{R}(Q)}\big]} = \overline{\mathcal{R}(P_{\mathcal{R}(Q)} - Q^*)}.$$

Replacing $Q$ with $I - Q$ in Lemma 2.3 gives

$$\mathcal{R}\big[m(I - Q)\big] = \mathcal{R}\big(|I - Q| + I - Q\big).$$

The equation above together with the second equation in (3.28) yields

$$\overline{\mathcal{R}(\widetilde{T_2}^*)} = \overline{\mathcal{R}\big[Q(|I - Q| + I - Q)\big]} = \overline{\mathcal{R}\big[Q|I - Q|\big]} = \overline{\mathcal{R}\big[Q(I - Q^*)\big]}$$
$$= \overline{\mathcal{R}\big[Q(I - P_{\mathcal{R}(Q)})\big]} = \overline{\mathcal{R}(P_{\mathcal{R}(Q)} - Q)}.$$

Hence, by Lemma 1.1 we see that $P_{\mathcal{R}(Q)} - Q$ has the polar decomposition if and only if $\overline{\mathcal{R}(T_1^*)}$ and $\overline{\mathcal{R}(\widetilde{T_2}^*)}$ are both orthogonally complemented in $H$. So, the desired conclusion can be derived from (i)$\Longleftrightarrow$(viii) in Theorem 3.7. □

*3.3. A way to construct semi-harmonious quasi-projection pairs*

It is worthwhile to find ways to construct semi-harmonious quasi-projection pairs. For this, we need a lemma as follow.

**Lemma 3.10.** *Suppose that $(P, Q)$ is a quasi-projection pair on $H$. Let $T_1$ and $T_2$ be defined by (3.1). Then*

$$\mathcal{R}(\lambda_1 T_1 + \lambda_2 T_2) = \mathcal{R}(T_1) + \mathcal{R}(T_2), \quad \forall \lambda_1, \lambda_2 \in \mathbb{C} \setminus \{0\}. \tag{3.29}$$

*Furthermore, the following statements are equivalent:*

(i) $(P, Q)$ *is semi-harmonious;*
(ii) $\forall \lambda_1, \lambda_2 \in \mathbb{C}\setminus\{0\}$, $\overline{\mathcal{R}(\lambda_1 T_1 + \lambda_2 T_2)}$ *is orthogonally complemented in $H$;*
(iii) $\exists \lambda_1, \lambda_2 \in \mathbb{C}\setminus\{0\}$, $\overline{\mathcal{R}(\lambda_1 T_1 + \lambda_2 T_2)}$ *is orthogonally complemented in $H$.*

*Proof.* Suppose that $\lambda_1 \ne 0$ and $\lambda_2 \ne 0$. Let

$$S = \lambda_1 T_1 + \lambda_2 T_2, \quad M_1 = \overline{\mathcal{R}(T_1)}, \quad M_2 = \overline{\mathcal{R}(T_2)}.$$

It is clear that

$$\mathcal{R}(S) \subseteq \mathcal{R}(\lambda_1 T_1) + \mathcal{R}(\lambda_2 T_2) = \mathcal{R}(T_1) + \mathcal{R}(T_2).$$

On the other hand, for every $x, y \in H$ we have

$$T_1 x + T_2 y = S\left[\frac{1}{\lambda_1}(I - Q)x + \frac{1}{\lambda_2}Qy\right] \in \mathcal{R}(S).$$

This shows the validity of (3.29).

Clearly, $M_1 \perp M_2$, so $M_1 + M_2$ is closed in $H$. Hence, by (3.29)
$$\overline{\mathcal{R}(S)} = M_1 + M_2. \qquad (3.30)$$
Therefore, $\overline{\mathcal{R}(S)}$ is invariant with respect to the parameters $\lambda_1, \lambda_2 \in \mathbb{C} \setminus \{0\}$, and thus (ii)$\Longleftrightarrow$(iii) follows.

Now, assume that $\overline{\mathcal{R}(S)}$ is orthogonally complemented in $H$. Then by (3.30) $H = M_1 + M_1^\perp$, where $M_1^\perp = M_2 + \overline{\mathcal{R}(S)}^\perp$. Hence, $M_1$ is orthogonally complemented in $H$. The same is true for $M_2$. Conversely, suppose that $M_1$ and $M_2$ are both orthogonally complemented in $H$. Then $P_{M_1} + P_{M_2}$ is a projection whose range is equal to $\overline{\mathcal{R}(S)}$ (see (3.30)). Therefore, $\overline{\mathcal{R}(S)}$ is orthogonally complemented in $H$. This shows that $\overline{\mathcal{R}(S)}$ is orthogonally complemented in $H$ if and only if $M_1$ and $M_2$ are both orthogonally complemented in $H$. Hence, the equivalence of (i) and (ii) follows immediately from (i)$\Longleftrightarrow$(ii) in Theorem 3.7. $\square$

Based on Lemma 3.10, we are now able to construct semi-harmonious quasi-projection pairs as follows.

**Theorem 3.11.** *For every quasi-projection pair $(P, Q)$ on $H$, let*
$$M = \overline{\mathcal{R}(T_1)} + \overline{\mathcal{R}(T_2)}, \qquad (3.31)$$
*where $T_1$ and $T_2$ are defined by (3.1). Then $(P|_M, Q|_M)$ is a semi-harmonious quasi-projection pair on $M$.*

*Proof.* Let $\lambda_1 = 1$ and $\lambda_2 = \pm 1$ in (3.29). From (3.30) and (3.1), we obtain
$$M = \overline{\mathcal{R}(T_1 + T_2)} = \overline{\mathcal{R}(T_1 - T_2)} = \overline{\mathcal{R}(P - Q)}, \qquad (3.32)$$
so $M$ is closed in $H$, hence it is a Hilbert $\mathfrak{A}$-module. Evidently,
$$\mathcal{R}(P|_M) = \overline{\mathcal{R}(T_1)} \subseteq M,$$
so $P|_M$ is a projection on $M$. By (3.1) and Lemma 2.2(ii), we have
$$\begin{aligned}
Q(P-Q) &= PQ(P-Q) + (I-P)Q(P-Q) \\
&= P[I - (I-Q)](P-Q) + T_2(P-Q) \\
&= T_1 - T_1(P-Q) + T_2(P-Q), \\
Q^*(P-Q) &= (2P-I)Q(2P-I) \cdot (P-Q) \\
&= PQ(2P-I)(P-Q) - (I-P)Q(2P-I)(P-Q) \\
&= P[I - (I-Q)](2P-I)(P-Q) - T_2(2P-I)(P-Q) \\
&= P(P-Q) - T_1(2P-I)(P-Q) - T_2(2P-I)(P-Q) \\
&= T_1 - T_1(2P-I)(P-Q) - T_2(2P-I)(P-Q).
\end{aligned}$$

From the above decompositions of $Q(P-Q)$ and $Q^*(P-Q)$, together with (3.31) and (3.32), we conclude that $Q|_M \in \mathcal{L}(M)$ such that $(Q|_M)^* = Q^*|_M$. Therefore, $(P|_M, Q|_M)$ is a quasi-projection pair on $M$ satisfying

$$P|_M(I_M - Q|_M) = T_1|_M, \quad (I_M - P|_M)Q|_M = T_2|_M.$$

From (3.31) and the first equation in (3.32), it is clear that

$$\overline{\mathcal{R}(T_1|_M + T_2|_M)} = \overline{\mathcal{R}\big[(T_1 + T_2)|_M\big]} = \overline{\mathcal{R}\big[(T_1 + T_2)^2\big]}.$$

Furthermore, $T_1 + T_2$ is self-adjoint (see Lemma 3.2), so from Lemma 1.2 we have

$$\overline{\mathcal{R}\big[(T_1+T_2)^2\big]} = \overline{\mathcal{R}\big[(T_1+T_2)(T_1+T_2)^*\big]} = \overline{\mathcal{R}(T_1+T_2)} = M.$$

This shows that $\overline{\mathcal{R}(T_1|_M + T_2|_M)} = M$, which obviously ensures the orthogonal complementarity of $\overline{\mathcal{R}(T_1|_M + T_2|_M)}$ in $M$. The desired conclusion is immediate from (i)$\Longleftrightarrow$(iii) in Lemma 3.10. $\square$

*Remark* 3.1. Given a quasi-projection pair $(P,Q)$ on $H$, let $T_i$ ($1 \le i \le 4$) and $M$ be defined by (3.1), (3.2) and (3.31), respectively. It is obvious that $\overline{\mathcal{R}(T_3|_M)} \subseteq \overline{\mathcal{R}(T_3)}$. On the other hand, by (3.9) we have

$$T_3 T_3^* = -T_3 T_4 = -PQ(I-P)QP,$$

which is combined with Lemma 1.2 and (3.31) to get

$$\overline{\mathcal{R}(T_3)} = \overline{\mathcal{R}(T_3 T_3^*)} = \overline{\mathcal{R}\big[(PQ(I-P)QP\big]} \subseteq \overline{\mathcal{R}\big[(PQ(I-P)Q\big]}$$
$$= \overline{\mathcal{R}(T_3 T_2)} = \overline{\mathcal{R}(T_3|_M T_2)} \subseteq \overline{\mathcal{R}(T_3|_M)}.$$

This shows that $\overline{\mathcal{R}(T_3|_M)} = \overline{\mathcal{R}(T_3)}$. Similarly, we have $\overline{\mathcal{R}(T_4|_M)} = \overline{\mathcal{R}(T_4)}$.

3.4. *Two types of non-semi harmonious quasi-projection pairs*

In the remainder of this section, we provide two types of quasi-projection pairs, both of which are not semi-harmonious.

**Example 3.1.** Suppose that $H_1$ and $H_2$ are Hilbert $C^*$-modules over a $C^*$-algebra, and $A \in \mathcal{L}(H_2, H_1)$ has no polar decomposition (the reader is referred to [10, Example 3.15] for such an operator). Let $Q \in \mathcal{L}(H_1 \oplus H_2)$ be the idempotent given by

$$Q = \begin{pmatrix} I_{H_1} & A \\ 0 & 0 \end{pmatrix}.$$

It is obvious that
$$P_{\mathcal{R}(Q)} = \begin{pmatrix} I_{H_1} & 0 \\ 0 & 0 \end{pmatrix}, \quad P_{\mathcal{R}(Q)} - Q = \begin{pmatrix} 0 & -A \\ 0 & 0 \end{pmatrix}.$$

In virtue of Lemma 1.1, we see that $P_{\mathcal{R}(Q)} - Q$ has no polar decomposition. So, by Theorem 3.9 $(m(Q), Q)$ is not semi-harmonious.

To construct another type of non semi-harmonious quasi-projection pairs, we focus our attention on the unital $C^*$-algebras. It is well-known that every $C^*$-algebra is a Hilbert module over itself. For the sake of completeness, we give a brief description. Suppose that $\mathfrak{B}$ is a $C^*$-algebra. Let
$$\langle x, y \rangle = x^* y, \quad \text{for } x, y \in \mathfrak{B}.$$

With the inner-product given as above, $\mathfrak{B}$ is a Hilbert $\mathfrak{B}$-module, and $\mathfrak{B}$ can be embedded into $\mathcal{L}(\mathfrak{B})$ via $b \to L_b$ [11, Section 3], where $L_b$ is defined by $L_b(x) = bx$ for $x \in \mathfrak{B}$. If $\mathfrak{B}$ has a unit, then it can be easily derived from (1.2) that $\mathcal{L}(\mathfrak{B}) = \{L_b : b \in \mathfrak{B}\}$. So, for every unital $C^*$-algebra $\mathfrak{B}$, we can identify $\mathfrak{B}$ with $\mathcal{L}(\mathfrak{B})$. To ease notation, we use $\mathcal{R}(b)$ and $\mathcal{N}(b)$ instead of $\mathcal{R}(L_b)$ and $\mathcal{N}(L_b)$, respectively. Given a unital $C^*$-algebra $\mathfrak{B}$, let $\mathfrak{A} = M_2(\mathfrak{B})$, which is a unital $C^*$-algebra. Hence, $\mathfrak{A}$ itself becomes a Hilbert $\mathfrak{A}$-module described as above.

*Remark* 3.2. Specifically, if $\mathfrak{B}$ is a unital commutative $C^*$-algebra, then up to $C^*$-isomorphic equivalence we may assume that $\mathfrak{B} = C(\Omega)$, the set consisting of all continuous functions on a compact Hausdorff space $\Omega$. Let $\mathfrak{A} = M_2(\mathfrak{B})$. Since $M_2(\mathbb{C})$ is nuclear and $\mathfrak{A} \cong M_2(\mathbb{C}) \otimes \mathfrak{B}$, the unique $C^*$-norm on $\mathfrak{A}$ is determined by
$$\|x\| = \max_{t \in \Omega} \|x(t)\| \quad (x \in \mathfrak{A}),$$
where $x_{ij} \in C(\Omega)$ $(1 \leq i, j \leq 2)$,
$$x = \begin{pmatrix} x_{11} & x_{12} \\ x_{21} & x_{22} \end{pmatrix}, \quad x(t) = \begin{pmatrix} x_{11}(t) & x_{12}(t) \\ x_{21}(t) & x_{22}(t) \end{pmatrix} \quad (t \in \Omega),$$
and for each $t \in \Omega$, $\|x(t)\|$ denotes the operator norm (also known as 2-norm) on $M_2(\mathbb{C})$.

**Theorem 3.12.** *Suppose that $\mathfrak{B}$ is a unital $C^*$-algebra and $A \in \mathfrak{B}_{sa}$ is such that $A^2 - A \geq 0$. Let $\mathfrak{A} = M_2(\mathfrak{B})$ and $\ell(A) \in \mathfrak{B}$ be given by (2.3), and let $P, Q \in \mathfrak{A}$ be defined by*
$$P = \begin{pmatrix} I_{\mathfrak{B}} & 0 \\ 0 & 0 \end{pmatrix}, \quad Q = \begin{pmatrix} A & -\ell(A) \\ \ell(A) & I_{\mathfrak{B}} - A \end{pmatrix}. \tag{3.33}$$

*Then $(P, Q)$ is semi-harmonious if and only if*

$$I_\mathfrak{B} \in \overline{\mathcal{R}(A - I_\mathfrak{B})} + \mathcal{N}(A - I_\mathfrak{B}), \qquad (3.34)$$

*where*

$$\mathcal{R}(A - I_\mathfrak{B}) = \{(A - I_\mathfrak{B})x : x \in \mathfrak{B}\},$$
$$\mathcal{N}(A - I_\mathfrak{B}) = \{x \in \mathfrak{B} : (A - I_\mathfrak{B})x = 0\}.$$

*Proof.* It is clear that $(P, Q)$ is a quasi-projection pair on $\mathfrak{A}$, and the operators $T_1$ and $T_2$ defined by (3.1) turn out to be

$$T_1 = \begin{pmatrix} I_\mathfrak{B} - A & \ell(A) \\ 0 & 0 \end{pmatrix}, \quad T_2 = \begin{pmatrix} 0 & 0 \\ \ell(A) & I_\mathfrak{B} - A \end{pmatrix}. \qquad (3.35)$$

Let

$$X = \mathcal{R}(A - I_\mathfrak{B}), \quad Y = \mathcal{N}(A - I_\mathfrak{B}).$$

Direct computations yield

$$(I_\mathfrak{B} - A)^2 + \ell^2(A) = (A - I_\mathfrak{B})(2A - I_\mathfrak{B}) = (2A - I_\mathfrak{B})(A - I_\mathfrak{B}).$$

By assumption $A^2 - A \ge 0$, which implies that $2A - I_\mathfrak{B}$ is invertible in $\mathfrak{B}$. So, from the above equations we obtain

$$\mathcal{R}\big[(I_\mathfrak{B} - A)^2 + \ell^2(A)\big] = X, \quad \mathcal{N}\big[(I_\mathfrak{B} - A)^2 + \ell^2(A)\big] = Y.$$

Therefore,

$$\overline{\mathcal{R}(T_1)} = \overline{\mathcal{R}(T_1 T_1^*)} = \left\{ \begin{pmatrix} x_{11} & x_{12} \\ 0 & 0 \end{pmatrix} : x_{11}, x_{12} \in \overline{X} \right\},$$

$$\mathcal{N}(T_1^*) = \mathcal{N}(T_1 T_1^*) = \left\{ \begin{pmatrix} x_{11} & x_{12} \\ x_{21} & x_{22} \end{pmatrix} : x_{11}, x_{12} \in Y; x_{21}, x_{22} \in \mathfrak{B} \right\}.$$

The equations above together with (1.2) indicate that

$$\mathfrak{A} = \overline{\mathcal{R}(T_1)} + \mathcal{N}(T_1^*) \iff \mathfrak{B} = \overline{X} + Y \iff I_\mathfrak{B} \in \overline{X} + Y.$$

This shows that $\overline{\mathcal{R}(T_1)}$ is orthogonally complemented in $\mathfrak{A}$ if and only if (3.34) is satisfied. From (3.35) and the derivations as above, it is easily seen that the same is true if $T_1$ is replaced by $T_2$. The desired conclusion follows from (i)$\iff$(ii) in Theorem 3.7. □

In view of Remark 3.2 and Theorem 3.12, it is easy to construct new type of non-semi harmonious quasi-projection pairs. We provide an example as follows.

**Example 3.2.** Let $\mathfrak{B} = C(\Omega)$ with $\Omega = [0,1]$ and $\mathfrak{A} = M_2(\mathfrak{B})$, which is regarded as the Hilbert module over itself. Let $P$ and $Q$ be defined by (3.33) such that the element $A$ in $\mathfrak{B}$ is given by

$$A(t) = \sec^2(t), \quad t \in [0,1].$$

Then $(A - I_\mathfrak{B})(t) = \tan^2(t)$ for any $t \in [0,1]$, whence

$$\mathcal{N}(A - I_\mathfrak{B}) = \{x \in C[0,1] : \tan^2(t)x(t) \equiv 0, t \in [0,1]\} = \{0\},$$

where $x(t) = 0$ for $t \neq 0$ comes from the observation that $\tan^2(t) \neq 0$ for all such $t$, and $x(0) = 0$ is derived from the continuity of $x$ at $t = 0$.

By definition, $\mathcal{R}(A - I_\mathfrak{B}) = \{(A - I_\mathfrak{B})y : y \in \mathfrak{B}\} = \{\tan^2(\cdot)y : y \in \mathfrak{B}\}$. In virtue of $\tan(0) = 0$, we have $x(0) = 0$ for every $x \in \overline{\mathcal{R}(A - I_\mathfrak{B})}$. Therefore, (3.34) fails to be true. Hence by Theorem 3.12, $(P,Q)$ is not semi-harmonious.

## 4. Harmonious quasi-projection pairs

The purpose of this section is to set up the general theory for harmonious quasi-projection pairs, which will play a crucial role in the block matrix representations for quasi-projection pairs [15, Section 7]. Inspired by [11, Definition 4.1] and [5, Theorem 2.4], we make a definition as follows.

**Definition 4.1.** A quasi-projection pair $(P,Q)$ on $H$ is referred to be harmonious if $H_5$ and $H_6$ defined by (3.6) are both orthogonally complemented in $H$.

*Remark* 4.1. Given a quasi-projection pair $(P,Q)$ on $H$, let $T_3$ and $T_4$ be defined by (3.2). We may combine Definition 4.1 with (3.6), (3.9) and Lemma 1.1 to conclude that

$$(P,Q) \text{ is harmonious} \Longleftrightarrow T_3 \text{ has the polar decomposition}$$
$$\Longleftrightarrow T_4 \text{ has the polar decomposition}.$$

As literally revealed, every harmonious quasi projection pair must be a semi harmonious quasi projection pair. Actually, a more in-depth characterization will be provided in Theorem 4.2. To achieve the desired characterization, we need a lemma as follows.

**Lemma 4.1.** *For every quasi-projection pair $(P,Q)$ on $H$, let $T_1$ and $H_5$ be defined by (3.1) and (3.6), respectively. Then the following statements are equivalent:*

(i) *$H_5$ is orthogonally complemented in $H$;*
(ii) *$\overline{\mathcal{R}(PQ)}$ and $\overline{\mathcal{R}(T_1)}$ are both orthogonally complemented in $H$.*

*In each case, $H_1$ and $H_2$ are both orthogonally complemented in $H$ such that*

$$P_{\overline{\mathcal{R}(PQ)}} = P_{H_1} + P_{H_5}, \quad P = P_{\overline{\mathcal{R}(PQ)}} + P_{H_2}, \quad P_{\overline{\mathcal{R}(T_1)}} = P_{H_2} + P_{H_5}. \quad (4.1)$$

*Proof.* Let $T_3, T_4, H_1$ and $H_2$ be defined by (3.2), (3.4) and (3.5), respectively.

(i)$\Longrightarrow$(ii). Assume that $H = H_5 + H_5^\perp$. From (3.6) and (3.9), we have

$$H_5^\perp = \overline{\mathcal{R}(T_3)}^\perp = \mathcal{N}(T_3^*) = \mathcal{N}(T_4).$$

So, a direct use of Lemma 3.5(ii) gives

$$H = H_1 + H_2 + H_5 + \mathcal{N}(P). \quad (4.2)$$

By the definitions of $H_1, H_2$ and $H_5$, it is clear that

$$\mathcal{N}(P) \perp H_i (i = 1, 2, 5), \quad H_5 \perp H_j (j = 1, 2).$$

Combining (3.15) and (3.5) yields $H_1 \perp H_2$. Hence, it can be concluded by (4.2) that $H_i (i = 1, 2, 5)$ are all orthogonally complemented in $H$ such that

$$P = P_{H_1} + P_{H_2} + P_{H_5}. \quad (4.3)$$

Rewrite (4.2) as
$$H = X_1 + X_2 = Y_1 + Y_2,$$

where $X_i, Y_i (i = 1, 2)$ are defined by

$$X_1 = H_2 + H_5, \quad X_2 = H_1 + \mathcal{N}(P), \quad Y_1 = H_1 + H_5, \quad Y_2 = H_2 + \mathcal{N}(P).$$

Note that $2P - I$ is surjective and $P(2P - I) = P$, so by Lemma 2.2(ii) we have $\mathcal{R}(PQ^*) = \mathcal{R}(PQ)$. It follows that

$$Y_1 \subseteq \overline{\mathcal{R}(PQ)}, \quad Y_2 \subseteq \mathcal{N}(QP) = \overline{\mathcal{R}(PQ^*)}^\perp = \overline{\mathcal{R}(PQ)}^\perp.$$

This shows that
$$H = \overline{\mathcal{R}(PQ)} + \overline{\mathcal{R}(PQ)}^\perp,$$

hence $\overline{\mathcal{R}(PQ)}$ is orthogonally complemented in $H$. Moreover,

$$P_{\overline{\mathcal{R}(PQ)}}H = P_{\overline{\mathcal{R}(PQ)}}(Y_1 + Y_2) = P_{\overline{\mathcal{R}(PQ)}}Y_1 = Y_1 = (P_{H_1} + P_{H_5})H,$$

so, the first equation in (4.1) is satisfied. In view of (4.3), the second equation in (4.1) follows. Similarly,

$$X_1 \subseteq \overline{\mathcal{R}(T_1)}, \quad X_2 \subseteq \mathcal{N}\big[(I-Q)P\big] = \overline{\mathcal{R}\big[P(I-Q)^*\big]}^\perp = \overline{\mathcal{R}(T_1)}^\perp.$$

Therefore, $\overline{\mathcal{R}(T_1)}$ is orthogonally complemented in $H$ such that the last equation in (4.1) is satisfied.

(ii)$\Longrightarrow$(i). From the first equation in (3.16), we have

$$H = \overline{\mathcal{R}(T_1)} + \mathcal{R}(I-P) + H_1. \tag{4.4}$$

By a replacement of $(P,Q)$ with $(P, I-Q)$, we obtain

$$H = \overline{\mathcal{R}(PQ)} + \mathcal{N}(P) + H_2. \tag{4.5}$$

It follows from (4.4) and (3.2) that

$$\begin{aligned}\mathcal{R}(PQ) \subseteq &\overline{\mathcal{R}\big[PQT_1\big]} + \mathcal{R}(T_3) + H_1 \\ =&\overline{\mathcal{R}\big[T_3(I-Q)\big]} + \mathcal{R}(T_3) + H_1 \subseteq H_5 + H_1.\end{aligned}$$

Therefore,

$$\overline{\mathcal{R}(PQ)} \subseteq H_5 + H_1.$$

This together with (4.5) yields the orthogonal decomposition $H = H_5 \dotplus Y$, where $Y$ is the right side of Lemma 3.5(ii). $\square$

**Theorem 4.2.** *For every quasi-projection pair $(P,Q)$ on $H$, $(P,Q)$ is harmonious if and only if $(P,Q)$ and $(P, I-Q)$ are both semi-harmonious.*

*Proof.* By a replacement of $(P,Q)$ with $(I-P, I-Q)$, it can be deduced from Lemma 4.1 that $H_6$ defined by (3.6) is orthogonally complemented in $H$ if and only if $\overline{\mathcal{R}\big[(I-P)(I-Q)\big]}$ and $\overline{\mathcal{R}(T_2)}$ are both orthogonally complemented in $H$. Thus, $(P,Q)$ is harmonious if and only if $\overline{\mathcal{R}(T_1)}$, $\overline{\mathcal{R}(T_2)}$, $\overline{\mathcal{R}(PQ)}$ and $\overline{\mathcal{R}\big[(I-P)(I-Q)\big]}$ are all orthogonally complemented in $H$. So, the desired conclusion follows from (i)$\Longleftrightarrow$(ii) in Theorem 3.7. $\square$

As a consequence of Theorem 4.2, we have the following result of symmetry.

**Corollary 4.3.** *For every quasi-projection pair $(P, Q)$ on $H$, let $\Sigma$ be defined by (2.1). If one element in $\Sigma$ is harmonious, then all the remaining elements in $\Sigma$ are harmonious.*

*Proof.* We may as well assume that $(P, Q)$ is harmonious. In this case, by Theorem 4.2 $(P, Q)$ and $(P, I - Q)$ are both semi-harmonious. Repeatedly using Corollary 3.8 twice indicates that the four pairs in (3.26) and these in

$$(P, I - Q), \quad (I - P, Q), \quad (P, I - Q^*) \quad \text{and} \quad (I - P, Q^*)$$

are all semi-harmonious. Hence, the conclusion is immediate from Theorem 4.2. □

*Remark* 4.2. Suppose that $(P, Q)$ is a harmonious quasi-projection pair. Let $T_1, T_2$ and $H_i(1 \le i \le 6)$ be defined by (3.1), (3.4), (3.5) and (3.6), respectively. By Definition 4.1, $H_5$ and $H_6$ are orthogonally complemented in $H$. So, by Lemma 4.1 $H_1$ and $H_2$ are orthogonally complemented in $H$ such that (4.1) is satisfied. Replacing $(P, Q)$ with $(I - P, Q)$, we see that $H_3$ and $H_4$ are also orthogonally complemented in $H$ such that

$$P_{\overline{\mathcal{R}(T_2)}} = P_{H_3} + P_{H_6}, \quad I - P = P_{\overline{\mathcal{R}(T_2)}} + P_{H_4},$$
$$P_{\overline{\mathcal{R}[(I-P)(I-Q)]}} = P_{H_4} + P_{H_6}.$$

We provide a featured result of the matched pair as follows.

**Theorem 4.4.** *For every matched pair $\bigl(m(Q), Q\bigr)$, it is harmonious if and only if it is semi-harmonious.*

*Proof.* Let $Q$ be an arbitrary idempotent on $H$. By Theorem 4.2, $\bigl(m(Q), Q\bigr)$ is harmonious if and only if both $\bigl(m(Q), Q\bigr)$ and $\bigl(m(Q), I - Q\bigr)$ are semi-harmonious. So, to get the conclusion, it is sufficient to prove that $\bigl(m(Q), I - Q\bigr)$ is always semi-harmonious. Let

$$T_1 = m(Q)Q, \quad T_2 = \bigl[I - m(Q)\bigr](I - Q).$$

From Lemmas 2.3, 1.4 and the invertibility of $|Q^*| + I$, we have

$$\overline{\mathcal{R}(T_1^*)} = \overline{\mathcal{R}\bigl[Q^*(|Q^*| + Q^*)\bigr]} = \overline{\mathcal{R}\bigl[Q^*(|Q^*| + I)\bigr]} = \mathcal{R}(Q^*),$$

which is orthogonally complemented in $H$. In virtue of the second equation stated in Lemma 2.4, it can be concluded that $\overline{\mathcal{R}(T_2^*)} = \mathcal{R}(I - Q^*)$ by a simple replacement of $Q$ with $I - Q$. Therefore, $\overline{\mathcal{R}(T_2^*)}$ is also orthogonally complemented in $H$. Hence, the desired conclusion can be derived immediately from (i)⟺(iii) in Theorem 3.7. □

As we know, every closed linear subspace of a Hilbert space is always orthogonally complemented. So, if $(P,Q)$ is a quasi-projection pair acting on a Hilbert space, then it is harmonious. It is notable that the same is true for some well-behaved Hilbert $C^*$-modules [12].

To construct new examples of harmonious and non-harmonious quasi-projection pairs, it is helpful to study furthermore the unital $C^*$-algebras, which are served as Hilbert $C^*$-modules.

**Theorem 4.5.** *Suppose that $\mathfrak{B}$ is a unital $C^*$-algebra and $A \in \mathfrak{B}_{sa}$ is such that $A^2 - A \geq 0$. Let $\mathfrak{A} = M_2(\mathfrak{B})$ and $\ell(A) \in \mathfrak{B}$ be given by (2.3), and let $P, Q \in \mathfrak{A}$ be defined by (3.33). Then $(P, Q)$ is harmonious if and only if*

$$I_\mathfrak{B} \in \overline{\mathcal{R}(A^2 - A)} + \mathcal{N}(A^2 - A), \tag{4.6}$$

*where*

$$\mathcal{R}(A^2 - A) = \{(A^2 - A)x : x \in \mathfrak{B}\},$$
$$\mathcal{N}(A^2 - A) = \{x \in \mathfrak{B} : (A^2 - A)x = 0\}.$$

*Proof.* Let $T_3$ and $T_4$ be defined by (3.2). Direct computations yield

$$T_3 = \begin{pmatrix} 0 & -\ell(A) \\ 0 & 0 \end{pmatrix}, \quad T_4 = \begin{pmatrix} 0 & 0 \\ \ell(A) & 0 \end{pmatrix}, \tag{4.7}$$

in which $\ell(A)$ is given by (2.3). Put

$$X = \mathcal{R}\big(\ell(A)\big), \quad Y = \mathcal{N}\big(\ell(A)\big).$$

Since $\ell(A) \geq 0$, by Lemmas 1.2 and 1.3 we have

$$\overline{X} = \overline{\mathcal{R}(\ell(A))} = \overline{\mathcal{R}(\ell^2(A))} = \overline{\mathcal{R}(A^2 - A)}, \quad Y = \mathcal{N}(\ell^2(A)) = \mathcal{N}(A^2 - A).$$

It follows from (4.7) that

$$H_5 = \overline{\mathcal{R}(T_3)} = \left\{ \begin{pmatrix} x_{11} & x_{12} \\ 0 & 0 \end{pmatrix} : x_{11}, x_{12} \in \overline{X} \right\},$$

$$H_5^\perp = \mathcal{N}(T_3^*) = \left\{ \begin{pmatrix} x_{11} & x_{12} \\ x_{21} & x_{22} \end{pmatrix} : x_{11}, x_{12} \in Y; x_{21}, x_{22} \in \mathfrak{B} \right\}.$$

Hence,
$$\mathfrak{A} = H_5 + H_5^\perp \iff \mathfrak{B} = \overline{X} + Y \iff I_\mathfrak{B} \in \overline{X} + Y.$$

This shows that $H_5$ is orthogonally complemented in $H$ if and only if (4.6) is satisfied. Similarly, it can also be concluded by (4.7) that $H_6$ is orthogonally complemented in $H$ if and only if (4.6) is satisfied. □

**Example 4.1.** Let $\mathfrak{B}$, $\mathfrak{A}$, $P$ and $Q$ be the same as in Example 3.2. As shown in Example 3.2, $(P,Q)$ fails to be semi-harmonious. Denote $I_{\mathfrak{A}}$ simply by $I$, and let $T_1, T_2, T_3, T_4$ and $M$ be defined by (3.1), (3.2) and (3.31), respectively. According to Theorem 3.11, $(P|_M, Q|_M)$ is semi-harmonious. In what follows, we prove that $(P|_M, Q|_M)$ is in fact a harmonious quasi-projection pair.

It should be aware that for a function $x$ in $C[0,1]$, the notation $\mathcal{R}(x)$ is usually employed for the range $\{x(t) : t \in [0,1]\}$ of $x$. To avoid such an ambiguity, it is helpful to use $\mathcal{R}(x(\cdot))$ to denote the set $\{xy : y \in C[0,1]\}$. Note that $\sec(t) \ne 0$ for every $t \in [0,1]$, so $\mathcal{R}[\sec(\cdot)\tan(\cdot)] = \mathcal{R}(\tan(\cdot))$. Therefore, by Lemma 1.2 we have

$$X := \overline{\mathcal{R}[\sec(\cdot)\tan(\cdot)]} = \overline{\mathcal{R}(\tan^2(\cdot))}. \tag{4.8}$$

Let $C, D \in \mathcal{L}(\mathfrak{B})$ be defined by

$$C = P(I-Q)P, \quad D = (I-P)Q(I-P).$$

It follows immediately from (2.1) and (2.2) that

$$C = C^* \quad \text{and} \quad D = D^*.$$

Since $C = T_1 P$, we have $\overline{\mathcal{R}(C)} \subseteq \overline{\mathcal{R}(T_1)}$. Utilizing (3.1) and (3.7) yields

$$T_1 T_1^* = P(I-Q)(2P-I)(I-Q)P = 2C^2 - C.$$

So, according to Lemma 1.2 we have

$$\overline{\mathcal{R}(T_1)} = \overline{\mathcal{R}(T_1 T_1^*)} = \overline{\mathcal{R}(2C^2 - C)} \subseteq \overline{\mathcal{R}(C)}.$$

Consequently,

$$\overline{\mathcal{R}(T_1)} = \overline{\mathcal{R}(2C^2 - C)} = \overline{\mathcal{R}(C)}. \tag{4.9}$$

Replacing $(P,Q)$ with $(I-P, I-Q)$ gives immediately

$$T_2 T_2^* = 2D^2 - D, \quad \overline{\mathcal{R}(T_2)} = \overline{\mathcal{R}(2D^2 - D)} = \overline{\mathcal{R}(D)}. \tag{4.10}$$

In view of (3.31), (4.9) and (4.10), we have $M = \overline{M_0}$, where

$$M_0 = \mathcal{R}(C) + \mathcal{R}(D).$$

A direct computation shows that for each $t \in [0,1]$,

$$C(t) = \begin{pmatrix} -\tan^2(t) & 0 \\ 0 & 0 \end{pmatrix}, \quad D(t) = \begin{pmatrix} 0 & 0 \\ 0 & -\tan^2(t) \end{pmatrix},$$

which means that for any $x \in \mathfrak{A}$ with $x(t) = \big(x_{ij}(t)\big)_{1\le i,j\le 2}$,

$$x \in M_0 \iff x_{ij} \in \mathcal{R}\big(\tan^2(\cdot)\big) \quad \text{for } i,j = 1,2.$$

Hence,
$$x \in M \iff x_{ij} \in X \quad \text{for } i,j = 1,2, \tag{4.11}$$

where $X$ is defined by (4.8).

From Remark 3.1 and (3.9), we see that the verification of the orthogonal decomposition
$$M = \overline{\mathcal{R}(T_4|_M)} + \mathcal{N}(T_4^*|_M)$$

can be simplified as
$$M = \overline{\mathcal{R}(T_4)} + \mathcal{N}(T_3) \cap M. \tag{4.12}$$

For every $x \in \mathfrak{A}$ with $x(t) = \big(x_{ij}(t)\big)_{1\le i,j\le 2}$, by (4.7) it is easy to verify that

$$x \in \mathcal{N}(T_3) \iff x_{21}(t) \equiv 0,\ x_{22}(t) \equiv 0.$$

So, it can be deduced from (4.11) that

$$x \in \mathcal{N}(T_3) \cap M \iff x_{21}(t) \equiv 0, x_{22}(t) \equiv 0; x_{11} \in X, x_{12} \in X. \tag{4.13}$$

Substituting $A(t) = \sec^2(t)$ into (4.7) yields

$$T_4(t) = \begin{pmatrix} 0 & 0 \\ \sec(t)\tan(t) & 0 \end{pmatrix}$$

for $t \in [0,1]$. Therefore, it can be concluded by (4.8) that

$$x \in \overline{\mathcal{R}(T_4)} \iff x_{11}(t) \equiv 0, x_{12}(t) \equiv 0; x_{21} \in X, x_{22} \in X. \tag{4.14}$$

Evidently, (4.12) can be derived immediately from (4.11), (4.13) and (4.14). This shows that $\overline{\mathcal{R}(T_4|_M)}$ is orthogonally complemented in $M$. The proof of the orthogonal complementarity of $\overline{\mathcal{R}(T_3|_M)}$ in $M$ is similar.

There exist semi-harmonious quasi-projection pairs, which fail to be harmonious. Our first example in this direction is as follows.

**Example 4.2.** Let $Q \in \mathcal{L}(H)$ be an idempotent such that $P_{\mathcal{R}(Q)} - Q$ has no polar decomposition. In this case, by Theorem 3.9 $\big(m(Q), Q\big)$ is not semi-harmonious, hence by Theorem 4.2 $\big(m(Q), I - Q\big)$ is not harmonious. However, from the proof of Theorem 4.4 we see that $\big(m(Q), I - Q\big)$ is semi-harmonious.

We provide another example as follows.

**Example 4.3.** Let $\mathfrak{B} = C(\Omega)$ with $\Omega = [-1, 0] \cup [1, 2]$ and $\mathfrak{A} = M_2(\mathfrak{B})$, which serves as the Hilbert module over itself. Let $P$ and $Q$ be defined by (3.33) such that the element $A$ in $\mathfrak{B}$ is specified as

$$A(t) = t, \quad \forall t \in \Omega.$$

Let $\ell$ be the function defined on $\Omega$ by $\ell(t) = \sqrt{t^2 - t}$. It is clear that

$$\ell(A)(t) = \ell(t), \quad \forall t \in \Omega,$$

where $\ell(A)$ is defined by (2.3). Therefore,

$$\ell(A) = \ell(\cdot) \in \mathfrak{B}. \tag{4.15}$$

Let $T_1, T_2, T_3, T_4$ and $M$ be defined by (3.1), (3.2) and (3.31), respectively. According to Theorem 3.11, $(P|_M, Q|_M)$ is semi-harmonious. In what follows, we prove that $(P|_M, Q|_M)$ is not harmonious.

Denote by $b = I_\mathfrak{B} - A$. Then

$$b(t) = 1 - t, \quad \forall t \in \Omega.$$

Let $Y, Z \subseteq \mathfrak{B}$ be defined by

$$Y = \overline{\mathcal{R}(b(\cdot))}, \quad Z = \overline{\mathcal{R}(\ell(\cdot))}$$

in the sense that for each $x \in \mathfrak{B}$,

$$\mathcal{R}(x(\cdot)) = \{xy : y \in \mathfrak{B}\} \subseteq \mathfrak{B}.$$

Observe that

$$Z = \overline{\mathcal{R}(\ell^2(\cdot))} \subseteq Y,$$

and $b = bI_\mathfrak{B} \in Y$ with $b(0) = 1$, whereas $z(0) = 0$ for any $z \in Z$. It follows that $Z \subsetneq Y$.

With the same convention, by (4.15) we have

$$\mathcal{N}[\ell(A)(\cdot)] = \mathcal{N}[\ell(\cdot)] = \{x \in \mathfrak{B} : \ell(t)x(t) \equiv 0, t \in \Omega\} = \{0\}. \tag{4.16}$$

Following the notations as in Example 4.1, we have $M = \overline{\mathcal{R}(C)} + \overline{\mathcal{R}(D)}$, where

$$C = P(I - Q)P = \begin{pmatrix} b & 0 \\ 0 & 0 \end{pmatrix}, \quad D = (I - P)Q(I - P) = \begin{pmatrix} 0 & 0 \\ 0 & b \end{pmatrix}.$$

So, for each $x \in \mathfrak{A}$ with $x(t) = \big(x_{ij}(t)\big)_{1 \le i,j \le 2}$,

$$x \in M \iff x_{ij} \in Y \quad \text{for } i,j = 1,2. \tag{4.17}$$

As shown in Example 4.1, $\overline{\mathcal{R}(T_4|_M)}$ is orthogonally complemented in $M$ if and only if (4.12) is satisfied. So, our aim here is to show that

$$\overline{\mathcal{R}(T_4)} + \mathcal{N}(T_3) \cap M \ne M. \tag{4.18}$$

For each $x \in \mathfrak{A}$ with the same expression as above, by (4.7), (4.15), (4.16) and (4.17), we conclude that

$$x \in \mathcal{N}(T_3) \cap M \iff x_{11}, x_{12} \in Y; x_{21}(t) \equiv 0, x_{22}(t) \equiv 0,$$
$$x \in \overline{\mathcal{R}(T_4)} \iff x_{11}(t) \equiv 0, x_{12}(t) \equiv 0; x_{21}, x_{22} \in Z.$$

Hence, (4.18) can be rephrased as $Z \ne Y$, which has already been verified.

## 5. Common similarity and unitary equivalence of operators

Inspired by [6, Section 3], in this section we consider mainly the common solvability of the following two operator equations

$$W(Q - P_{H_1})W^{-1} = I - P - P_{H_4}, \tag{5.1}$$
$$W(P - P_{H_1})W^{-1} = I - Q - P_{H_4}, \tag{5.2}$$

where $(P, Q)$ is a known quasi-projection pair on $H$ such that $H_1$ and $H_4$ defined by (3.4) are both orthogonally complemented in $H$, while $W \in \mathcal{L}(H)$ is an unknown invertible operator.

When $(P, Q)$ is semi-harmonious, a concrete invertible operator $W$ will be constructed in Theorem 5.3 such that (5.1) and (5.2) are both satisfied. To achieve this main result of this section, we need a couple of lemmas.

**Lemma 5.1.** *Suppose that $(P, Q)$ is a semi-harmonious quasi-projection pair on $H$. Let $T_1$, $T_2$, $H_1$ and $H_4$ be defined by (3.1) and (3.4) respectively, and let $T_i = V_i|T_i|$ $(i = 1, 2)$ be the polar decompositions. Then*

$$P_{H_1}(P - P_{H_1}) = (I - Q - P_{H_4})P_{H_1} = 0, \tag{5.3}$$
$$P_{H_4}(P - P_{H_1}) = (I - Q - P_{H_4})P_{H_4} = 0. \tag{5.4}$$

*Furthermore, for every $\lambda_1, \lambda_2 \in \mathbb{C}$ the following statements are equivalent:*

(i) $(\lambda_1 V_1 - \lambda_2 V_2)(P - P_{H_1}) = (I - Q - P_{H_4})(\lambda_1 V_1 - \lambda_2 V_2);$

(ii) $\lambda_i V_i T_j = \lambda_j T_i V_j$ for $i,j = 1,2$.

*Proof.* By Theorem 3.7 $P_{H_1}$ and $P_{H_4}$ (henceforth denoted simply by $P_1$ and $P_4$) are existent. Evidently, (5.3) and (5.4) can be derived directly from (3.4). So, it remains to prove the equivalence of (i) and (ii).

From (3.16) and (3.25), we have

$$P - P_1 = P_{\overline{\mathcal{R}(T_1)}} = V_1 V_1^*, \tag{5.5}$$
$$(I - P) - P_4 = P_{\overline{\mathcal{R}(T_2)}} = V_2 V_2^*. \tag{5.6}$$

Observe that

$$\mathcal{R}(V_1) = \overline{\mathcal{R}(T_1)} \subseteq \mathcal{R}(P), \quad \mathcal{R}(V_2) = \overline{\mathcal{R}(T_2)} \subseteq \mathcal{N}(P),$$

so by (5.5) and (5.6) we can get

$$V_2^* V_1 = V_2^*(V_2 V_2^* V_1 V_1^*) V_1 = V_2^* \cdot 0 \cdot V_1 = 0,$$
$$P_1 V_1 = P_1 V_2 = 0, \quad P_4 V_1 = P_4 V_2 = 0,$$

which yield

$$V_1^* V_2 = 0, \quad V_1^* P_1 = V_2^* P_1 = 0, \quad V_1^* P_4 = V_2^* P_4 = 0 \tag{5.7}$$

by taking $*$-operation. Also, from the definitions of $P_1$, $P_4$, $T_1$ and $T_2$ it is clear that

$$T_1 P_1 = T_2 P_1 = T_1 P_4 = T_2 P_4 = T_1 Q = T_2(I - Q) = 0,$$

which implies that

$$V_1 P_1 = V_2 P_1 = V_1 P_4 = V_2 P_4 = V_1 Q = V_2(I - Q) = 0, \tag{5.8}$$

since it is true that $\mathcal{N}(V_i) = \mathcal{N}(T_i)$ for $i = 1,2$. It follows that

$$V_1 = P V_1 (I - Q), \quad V_2 = (I - P) V_2 Q. \tag{5.9}$$

Given $\lambda_1, \lambda_2 \in \mathbb{C}$, let

$$A = (\lambda_1 V_1 - \lambda_2 V_2)(P - P_1), \quad B = (I - Q - P_4)(\lambda_1 V_1 - \lambda_2 V_2).$$

Then by (5.5)–(5.9), we have

$$\begin{aligned}
A &= (\lambda_1 V_1 - \lambda_2 V_2) P = \lambda_1 P V_1 (I - Q) P - \lambda_2 (I - P) V_2 Q P, \\
B &= \big[(P + V_2 V_2^*) - Q\big](\lambda_1 V_1 - \lambda_2 V_2) = (P - Q)(\lambda_1 V_1 - \lambda_2 V_2) - \lambda_2 V_2 \\
&= \lambda_1 V_1 - Q(\lambda_1 V_1 - \lambda_2 V_2) - \lambda_2 V_2 = (I - Q)(\lambda_1 V_1 - \lambda_2 V_2) \\
&= (I - Q)\big[\lambda_1 P V_1 (I - Q) - \lambda_2 (I - P) V_2 Q\big].
\end{aligned}$$

Consequently,

$$P(A-B)Q = \lambda_1 V_1(I-Q)PQ + \lambda_2 P(I-Q)V_2 = -\lambda_1 V_1 T_2 + \lambda_2 T_1 V_2.$$

Therefore,
$$P(A-B)Q = 0 \iff \lambda_1 V_1 T_2 = \lambda_2 T_1 V_2.$$

Similarly, it can be verified that

$$P(A-B)(I-Q) = 0 \iff \lambda_1 V_1 T_1 = \lambda_1 T_1 V_1,$$
$$(I-P)(A-B)Q = 0 \iff \lambda_2 V_2 T_2 = \lambda_2 T_2 V_2,$$
$$(I-P)(A-B)(I-Q) = 0 \iff \lambda_2 V_2 T_1 = \lambda_1 T_2 V_1.$$

This shows the equivalence of items (i) and (ii). $\square$

**Lemma 5.2.** *Under the notations and conditions of Lemma 5.1, the following statements are equivalent:*

(i) $(V_1 + V_2)(T_1 + T_2) = (T_1 + T_2)(V_1 + V_2)$;
(ii) $V_i T_j = T_i V_j$ for $i, j = 1, 2$;
(iii) $|T_i| \cdot |T_j^*|^2 = |T_i|^2 \cdot |T_j^*|$ for $i, j = 1, 2$.

*Proof.* (i)$\iff$(ii). Let $A = (V_1 + V_2)(T_1 + T_2)$ and $B = (T_1 + T_2)(V_1 + V_2)$. By (5.9) and (3.1), we have

$$(I-P)V_1 = PV_2 = V_1 Q = V_2(I-Q) = 0, \tag{5.10}$$
$$(I-P)T_1 = PT_2 = T_1 Q = T_2(I-Q) = 0. \tag{5.11}$$

So, the equation $V_1 T_2 = T_1 V_2$ turns out to be $PAQ = PBQ$. Similarly, the other three equations in (ii) can be rephrased as $XAY = XBY$ for all the cases that $X \in \{P, I-P\}$ and $Y \in \{Q, I-Q\}$ satisfying $(X, Y) \neq (P, Q)$. Hence, the equivalence of (i) and (ii) is derived.

(ii)$\iff$(iii). In virtue of $\mathcal{R}(V_j) = \overline{\mathcal{R}(|T_j^*|)}$ and $\mathcal{N}(V_i) = \mathcal{N}(|T_i|)$ for $i, j = 1, 2$, we have

$$V_i T_j - T_i V_j = 0 \iff V_i(|T_j^*| - |T_i|)V_j = 0$$
$$\iff V_i(|T_j^*| - |T_i|)|T_j^*| = 0$$
$$\iff |T_i|(|T_j^*| - |T_i|)|T_j^*| = 0.$$

This shows the equivalence of (ii) and (iii). $\square$

*Remark* 5.1. By Lemma 3.3, we know that Lemma 5.2(iii) is always true for every quasi-projection pair $(P,Q)$. However, when $(P,Q)$ is not semi-harmonious, the operators $T_1$ and $T_2$ defined by (3.1) may have no polar decompositions. In such case, items (i) and (ii) in Lemma 5.2 are meaningless.

Recall that an element $A \in \mathcal{L}(H)$ is called an EP operator if $\mathcal{R}(A) = \mathcal{R}(A^*)$. If $A$ is Moore-Penrose invertible (equivalently, $\mathcal{R}(A)$ is closed in $H$ [19, Theorem 2.2]), then this can be characterized as $AA^\dagger = A^\dagger A$, since $AA^\dagger$ and $A^\dagger A$ are projections whose ranges are equal to $\mathcal{R}(A)$ and $\mathcal{R}(A^*)$ respectively, where $A^\dagger$ denotes the Moore-Penrose inverse of $A$. It is notable that an EP operator needs not to be self-adjoint. When $(P,Q)$ is a semi-harmonious quasi-projection pair, let $V_1$ and $V_2$ be defined as above. In the next theorem, we will show that $V_1 - V_2$ is an EP operator, which plays a crucial role in the verification of the invertibility of the operator $W$ constructed as (5.12).

**Theorem 5.3.** *Suppose that $(P,Q)$ is a semi-harmonious quasi-projection pair on $H$. Let $T_1$, $T_2$, $H_1$ and $H_4$ be defined by (3.1) and (3.4) respectively, and let $T_i = V_i|T_i|$ $(i=1,2)$ be the polar decompositions. Then $V_1 + V_2$ is self-adjoint, and the operator $W \in \mathcal{L}(H)$ given by*

$$W = \lambda_1(V_1 - V_2) + \lambda_2 P_{H_1} + \lambda_3 P_{H_4} \tag{5.12}$$

*is a common solution of (5.1) and (5.2), where $\lambda_i(i=1,2,3)$ are arbitrary in $\mathbb{C} \setminus \{0\}$.*

*Proof.* By Theorem 3.7, $P_{H_1}$ and $P_{H_4}$ (simplified as $P_1$ and $P_4$) are existent. From Lemmas 5.1, 5.2 and 3.3, we have

$$(\lambda_1 V_1 - \lambda_1 V_2)(P - P_1) = (I - Q - P_4)(\lambda_1 V_1 - \lambda_1 V_2),$$

which is combined with (5.3), (5.4) and (5.12) to get

$$W(P - P_1) = (I - Q - P_4)W.$$

Observe that $V_2^* V_1 = 0$ (see the proof of Lemma 5.1), so from (5.6), (5.12), (5.7), (3.15), (5.10) and (5.8) we can obtain

$$(I - P - P_4)W = V_2 V_2^* W = -\lambda_1 V_2,$$
$$W(Q - P_1) = \lambda_1(V_1 - V_2)(Q - P_1) = -\lambda_1 V_2.$$

Hence,
$$W(Q - P_1) = (I - P - P_4)W.$$

So, it remains only to prove that $W$ is invertible and $V_1 + V_2$ is self-adjoint.

Firstly, we prove that $V_1 - V_2$ is an EP operator. By (5.11), we have
$$T_1 = P(T_1 \pm T_2), \quad T_2 = (I - P)(T_2 \pm T_1),$$
$$T_1^* = (I - Q^*)(T_1^* \pm T_2^*), \quad T_2^* = Q^*(T_2^* \pm T_1^*),$$

which lead obviously to
$$\mathcal{N}(T_1) \cap \mathcal{N}(T_2) = \mathcal{N}(T_1 \pm T_2), \quad \mathcal{N}(T_1^* \pm T_2^*) = \mathcal{N}(T_1^*) \cap \mathcal{N}(T_2^*).$$

As $T_1 + T_2$ is self-adjoint (see Lemma 3.2), this gives
$$\mathcal{N}(T_1) \cap \mathcal{N}(T_2) = \mathcal{N}(T_1 + T_2) = \mathcal{N}\big[(T_1 + T_2)^*\big] = \mathcal{N}(T_1^*) \cap \mathcal{N}(T_2^*).$$

Hence
$$\mathcal{N}(V_1) \cap \mathcal{N}(V_2) = \mathcal{N}(V_1^*) \cap \mathcal{N}(V_2^*).$$

The same method employed as above shows that
$$\mathcal{N}(V_1 \pm V_2) = \mathcal{N}(V_1) \cap \mathcal{N}(V_2) = \mathcal{N}(V_1^*) \cap \mathcal{N}(V_2^*) = \mathcal{N}\big[(V_1 \pm V_2)^*\big]. \quad (5.13)$$

Following the line in the derivation of (3.29), we have
$$\mathcal{R}(\lambda_1 V_1 + \lambda_2 V_2) = \mathcal{R}(V_1) + \mathcal{R}(V_2), \quad \forall \lambda_1, \lambda_2 \in \mathbb{C} \setminus \{0\}.$$

Combining the observation $V_1 V_1^* \perp V_2 V_2^*$ with (5.5) and (5.6) yields
$$\mathcal{R}(V_1) + \mathcal{R}(V_2) = \mathcal{R}(V_1 V_1^*) + \mathcal{R}(V_2 V_2^*) = \mathcal{R}(V_1 V_1^* + V_2 V_2^*) = \mathcal{R}(\widetilde{P}),$$

where $\widetilde{P}$ is a projection on $H$ given by
$$\widetilde{P} = I - P_1 - P_4. \quad (5.14)$$

Specifically, if we let $\lambda_1 = 1$ and $\lambda_2 = -1$, then it gives
$$\mathcal{R}(V_1 - V_2) = \mathcal{R}(\widetilde{P}),$$

and thus $\mathcal{R}(V_1 - V_2)$ is closed in $H$, so $V_1 - V_2$ is Moore-Penrose invertible such that
$$(V_1 - V_2)(V_1 - V_2)^\dagger = \widetilde{P}. \quad (5.15)$$

It follows from (5.15) and (5.13) that
$$\mathcal{N}(\widetilde{P}) = \mathcal{N}\big[(V_1 - V_2)^\dagger\big] = \mathcal{N}\big[(V_1 - V_2)^*\big] = \mathcal{N}(V_1 - V_2)$$
$$= \mathcal{N}\left[(V_1 - V_2)^\dagger (V_1 - V_2)\right],$$

which implies that
$$(V_1 - V_2)^\dagger (V_1 - V_2) = \widetilde{P}. \tag{5.16}$$

Hence, by (5.15) and (5.16) $V_1 - V_2$ is an EP operator.

Secondly, we prove that $W^{-1} = \widetilde{W}$, where
$$\widetilde{W} = \frac{1}{\lambda_1}(V_1 - V_2)^\dagger + \frac{1}{\lambda_2}P_1 + \frac{1}{\lambda_3}P_4.$$

Indeed, by (5.14) we have $\widetilde{P}P_1 = P_1\widetilde{P} = 0$, which is combined with (5.15) and (5.16) to get
$$(V_1 - V_2)^\dagger P_1 = (V_1 - V_2)^\dagger \left[(V_1 - V_2)(V_1 - V_2)^\dagger\right] P_1 = (V_1 - V_2)^\dagger (\widetilde{P}P_1) = 0,$$
$$P_1(V_1 - V_2)^\dagger = P_1 \left[(V_1 - V_2)^\dagger(V_1 - V_2)\right](V_1 - V_2)^\dagger = (P_1\widetilde{P})(V_1 - V_2)^\dagger = 0.$$

The same is true if $P_1$ is replaced by $P_4$. That is,
$$(V_1 - V_2)^\dagger P_4 = P_4(V_1 - V_2)^\dagger = 0.$$

With the derivations as above and the various equations given in the proof of Lemma 5.1, it is easily seen that $W\widetilde{W} = \widetilde{W}W = I$.

Finally, we prove that $V_1 + V_2$ is self-adjoint. To this end, we show that
$$(V_1 + V_2)^* \widetilde{P} = (V_1 + V_2)\widetilde{P}, \quad (V_1 + V_2)^*(I - \widetilde{P}) = (V_1 + V_2)(I - \widetilde{P}), \tag{5.17}$$

where $\widetilde{P}$ is given by (5.14). By (5.15) and (5.16), we have
$$\mathcal{N}(\widetilde{P}) = \mathcal{N}(V_1 - V_2) = \mathcal{N}\left[(V_1 - V_2)^*\right],$$

which is combined with (5.13) to get
$$\mathcal{R}(I - \widetilde{P}) = \mathcal{N}(\widetilde{P}) = \mathcal{N}(V_1 + V_2) = \mathcal{N}\left[(V_1 + V_2)\right]^*.$$

Therefore, the second equation in (5.17) is satisfied.

Note that $\mathcal{R}(\widetilde{P}) = \mathcal{R}(V_1) + \mathcal{R}(V_2) = \overline{\mathcal{R}(T_1 + T_2)}$ (see (3.30)), the first equation in (5.17) can therefore be changed into
$$(V_1 + V_2)^*(T_1 + T_2) = (V_1 + V_2)(T_1 + T_2),$$

which can be rephrased furthermore as
$$(V_1 + V_2)^*(T_1 + T_2) = (T_1 + T_2)(V_1 + V_2) \tag{5.18}$$

by using (i)⟺(iii) in Lemma 5.2 and Lemma 3.3. As observed before, we have $V_2^* V_1 = 0$, which gives $V_2^* T_1 = 0$ since $\mathcal{R}(T_1)$ is contained in $\mathcal{R}(V_1)$. Similarly, $V_1^* T_2 = 0$. As a result,

$$(V_1 + V_2)^*(T_1 + T_2) = V_1^* T_1 + V_2^* T_2 = |T_1| + |T_2|,$$

which leads clearly to (5.18) by taking $*$-operation, since both $|T_1| + |T_2|$ and $T_1 + T_2$ are self-adjoint. □

Let $Q \in \mathcal{L}(H)$ be an idempotent. It is well-known that $\|Q\| > 1$ whenever $Q$ is not a projection. Under the condition of Theorem 5.3 and suppose furthermore that $Q$ is not a projection, $I - Q - P_4$ is an idempotent which is not a projection, so

$$\|I - Q - P_{H_4}\| > 1 \geq \|P - P_{H_1}\| = \|U(P - P_{H_1})U^*\|$$

for any unitary operator $U$. Therefore, the operator $W$ constructed in Theorem 5.3 is merely an invertible operator, which generally fails to be a unitary operator. Employing (5.1) and (5.2) yields

$$W(P - P_{H_1})(Q - P_{H_1})W^{-1} = (I - Q - P_{H_4})(I - P - P_{H_4}). \qquad (5.19)$$

Our next theorem indicates that the operator $W$ in the above modified equation can in fact be replaced by a unitary operator.

**Theorem 5.4.** *Under the condition and the notations of Theorem 5.3, we have*

$$U(P - P_{H_1})(Q - P_{H_1})U^* = (I - Q - P_{H_4})(I - P - P_{H_4}), \qquad (5.20)$$

*where $U$ is a unitary operator defined by*

$$U = \lambda_1 \left[V_1(2P - I) - V_2\right] + \lambda_2 P_1 + \lambda_3 P_4, \qquad (5.21)$$

*in which $\lambda_i (i = 1, 2, 3)$ are arbitrary in the unit circle of $\mathbb{C}$.*

*Proof.* As before, denote $P_{H_1}$ and $P_{H_4}$ simply by $P_1$ and $P_4$, respectively. Suppose that $|\lambda_i| = 1$ for $i = 1, 2, 3$. First, we prove that the operator $U$ defined by (5.21) is a unitary. Let $\widetilde{T_1}$ be defined by (3.3) and let

$$\widetilde{V_1} = V_1(2P - I).$$

Since $2P - I$ is a symmetry, we have

$$\widetilde{V_1}\widetilde{V_1}^* = V_1(2P-I)(2P-I)V_1^* = V_1V_1^* = P_{\overline{\mathcal{R}(T_1)}} = P_{\overline{\mathcal{R}(\widetilde{T_1})}},$$

$$\left|\widetilde{T_1}^*\right|^2 = T_1(2P-I)(2P-I)T_1^* = T_1T_1^* = |T_1^*|^2,$$

and thus $\left|\widetilde{T_1}^*\right| = |T_1^*|$. Also, it is clear that

$$\widetilde{T_1}^* = (2P-I)T_1^* = (2P-I)V_1^*|T_1^*| = \widetilde{V_1}^*\left|\widetilde{T_1}^*\right|.$$

So, $\widetilde{T_1}^* = \widetilde{V_1}^*\left|\widetilde{T_1}^*\right|$ is the polar decomposition. In addition, by (3.7) we have

$$T_2\widetilde{T_1}^* = (I-P)Q(I-Q)P = 0,$$

which gives $V_2\widetilde{V_1}^* = 0$, since $\overline{\mathcal{R}(\widetilde{T_1}^*)} = \mathcal{R}(\widetilde{V_1}^*)$ and $\mathcal{N}(T_2) = \mathcal{N}(V_2)$. Hence, $V_2^*V_2 \perp \widetilde{V_1}^*\widetilde{V_1}$. Moreover, by (3.25), (5.5) and (5.6) we have

$$\widetilde{V_1}^*\widetilde{V_1} = P_{\overline{\mathcal{R}(\widetilde{T_1}^*)}} = P_{\mathcal{N}(Q)} - P_4,$$

$$\widetilde{V_1}\widetilde{V_1}^* + V_2V_2^* = V_1V_1^* + V_2V_2^* = I - P_1 - P_4.$$

Consequently,

$$UU^* = \widetilde{V_1}\widetilde{V_1}^* + V_2V_2^* + P_1 + P_4 = I,$$
$$U^*U = \widetilde{V_1}^*\widetilde{V_1} + V_2^*V_2 + P_1 + P_4 = P_{\mathcal{N}(Q)} + V_2^*V_2 + P_1.$$

So, $U^*U$ is a summation of three mutually orthogonal projections. Thus, in order to get $U^*U = I$ it is sufficient to verify the triviality of the common null space of these projections. It is routine to make such a verification. For, if $x \in H$ satisfies

$$T_2x = 0, \quad P_1x = 0, \quad P_{\mathcal{N}(Q)}x = 0,$$

then $Qx = PQx \in H_1$, so $Qx = P_1Qx = P_1x = 0$, and therefore $x = P_{\mathcal{N}(Q)}x = 0$.

Next, we check the validity of

$$X(P-P_1)(Q-P_1) = (I-Q-P_4)(I-P-P_4)X. \qquad (5.22)$$

in the case that $X = U$. From (5.19) we see that $W$ is a solution of the above equation. If $X = V_1(P-I)$, then both sides of (5.22) are equal to zero. So, $V_1(P-I)$ is another solution. This finishes the verification of (5.22) for $X = U$, since by (5.12) and (5.21) we have $U = W + 2\lambda_1 V_1(P-I)$. □

Some new phenomena may happen if a quasi-projection pair fails to be semi-harmonious. In fact, we have a theorem as follows.

**Theorem 5.5.** *There exists a quasi-projection pair $(P, Q)$ on certain Hilbert $C^*$-module $H$ such that $H_1$ and $H_4$ defined by (3.4) satisfy $H_1 = \{0\}$ and $H_4 = \{0\}$, whereas (5.22) has no solution of invertible operator.*

*Proof.* We use the framework described in Theorem 3.12 and make a slight modification of that in Example 3.2. Let $\mathfrak{B} = C(\Omega)$ with $\Omega = [0, 1]$ and $\mathfrak{J}$ be the ideal of $\mathfrak{B}$ given by

$$\mathfrak{J} = \{f \in \mathfrak{B} : f(0) = 0\}.$$

Let $\mathfrak{A}$ be the unital $C^*$-subalgebra of $M_2(\mathfrak{B})$ defined by

$$\mathfrak{A} = \left\{ \begin{pmatrix} b_{11} & b_{12} \\ b_{21} & b_{22} \end{pmatrix} : b_{11}, b_{22} \in \mathfrak{B}; b_{12}, b_{21} \in \mathfrak{J} \right\}.$$

As before, $\mathfrak{A}$ serves as a Hilbert module over itself. Let $P, Q \in \mathfrak{A}$ be the same as that in Example 3.2. Namely,

$$P(t) = \begin{pmatrix} 1 & 0 \\ 0 & 0 \end{pmatrix}, \quad Q(t) = \begin{pmatrix} \sec^2(t) & -\sec(t)\tan(t) \\ \sec(t)\tan(t) & -\tan^2(t) \end{pmatrix}$$

for each $t \in \Omega$. Suppose that $x \in H_4$ is determined by $x(t) = \big(x_{ij}(t)\big)_{1 \le i,j \le 2}$ for each $t \in \Omega$. Then

$$P(t)x(t) \equiv 0, \quad Q(t)x(t) \equiv 0 \quad (t \in \Omega),$$

which clearly lead to $x_{i,j}(t) \equiv 0$ for all $i, j = 1, 2$, since $\sec(t) \ne 0$ for every $t \in \Omega$ and $\tan(t) = 0$ only if $t = 0$, and all the functions considered are continuous on $\Omega$. This shows that $H_4 = \{0\}$. Similarly, $H_1 = \mathcal{N}(I - P) \cap \mathcal{N}(I - Q) = \{0\}$. Therefore, (5.22) is simplified as

$$WPQ = (I - Q)(I - P)W. \tag{5.23}$$

Now, suppose that $W \in \mathfrak{A}$ is determined by $W(t) = \big(W_{ij}(t)\big)_{1 \le i,j \le 2}$ and is a solution of (5.23). Substituting $t = 0$ into (5.23) yields

$$\begin{pmatrix} W_{11}(0) & W_{12}(0) \\ W_{21}(0) & W_{22}(0) \end{pmatrix} \begin{pmatrix} 1 & 0 \\ 0 & 0 \end{pmatrix} = \begin{pmatrix} 0 & 0 \\ 0 & 1 \end{pmatrix} \begin{pmatrix} W_{11}(0) & W_{12}(0) \\ W_{21}(0) & W_{22}(0) \end{pmatrix}.$$

Hence, $W_{11}(0) = W_{22}(0) = 0$. In addition, since $W_{12}$ and $W_{21}$ are taken in $\mathfrak{J}$ which requires $W_{12}(0) = W_{21}(0) = 0$. As a result, $W(0)$ is the zero matrix in $M_2(\mathbb{C})$, which clearly implies the non-invertibility of $W$ in $\mathfrak{A}$. □

As an immediate consequence of Theorem 5.5, we have the following corollary.

**Corollary 5.6.** *There exists a quasi-projection pair $(P,Q)$ on certain Hilbert $C^*$-module $H$ such that $H_1$ and $H_4$ defined by (3.4) satisfy $H_1 = \{0\}$ and $H_4 = \{0\}$, whereas (5.1) and (5.2) have no common solution of invertible operator.*

## 6. A norm equation associated with the Friedrichs angle

Let $M$ and $N$ be closed subspaces of a Hilbert space. The Friedrichs angle [4], denoted by $\alpha(M,N)$, is the unique angle in $[0, \frac{\pi}{2}]$ whose cosine is equal to $c(M,N)$, where

$$c(M,N) = \sup\{|\langle x,y\rangle| : x \in \widetilde{M}, y \in \widetilde{N}, \|x\| \le 1, \|y\| \le 1\},$$

in which $\widetilde{M} = M \cap (M \cap N)^\perp$ and $\widetilde{N} = N \cap (M \cap N)^\perp$. It is known [3, Lemma 10] that

$$c(M,N) = \|P_M P_N - P_{M \cap N}\|, \tag{6.1}$$

which is invariant in the sense that

$$c(M,N) = c(M^\perp, N^\perp). \tag{6.2}$$

Inspired by (6.1) and (6.2), the following lemma is obtained in [11] for two projections on a general Hilbert $C^*$-module.

**Lemma 6.1.** *[11, Theorem 5.12] Let $P, Q \in \mathcal{L}(H)$ be projections such that $H_1$ and $H_4$ defined by (3.4) are orthogonally complemented in $H$. Then*

$$\|PQ - P_{H_1}\| = \|(I-P)(I-Q) - P_{H_4}\|. \tag{6.3}$$

The characteristics of the above lemma is mainly reflected in its weakest condition, since for each $i \in \{1,4\}$, $P_{H_i}$ exists if and only if $H_i$ is orthogonally complementary in $H$. It is interesting to investigate the validity of (6.3) in the new setting of quasi-projection pairs on a general Hilbert $C^*$-module. To give the positive answer under the above weakest condition, we need a couple of lemmas, which are provided as follows.

**Lemma 6.2.** *Equation (6.3) is valid for every semi-harmonious quasi-projection pair $(P,Q)$ on $H$, where $H_1$ and $H_4$ are defined by (3.4).*

*Proof.* Let $(P, Q)$ be an arbitrary quasi-projection pair on $H$. From the second equations in (3.8), (3.3) and (3.1), we obtain

$$[Q(I-P)]^* = \widetilde{T_2} = (I-P)Q(I-2P).$$

Replacing $Q$ with $I - Q$ gives

$$[(I-Q)(I-P)]^* = (I-P)(I-Q)(I-2P). \qquad (6.4)$$

Now, suppose that $H_1$ and $H_4$ defined by (3.4) are both orthogonally complemented in $H$. Let $P_{H_1}$ and $P_{H_4}$ be denoted simply by $P_1$ and $P_4$, respectively. By (6.4) and the observation $P_4^* = P_4 = P_4(I - 2P)$, we have

$$[(I-Q)(I-P) - P_4]^* = [(I-P)(I-Q) - P_4](I-2P).$$

Since $I - 2P$ is a unitary, this shows that

$$\|[(I-Q)(I-P) - P_4]^*\| = \|(I-P)(I-Q) - P_4\|. \qquad (6.5)$$

Suppose furthermore that $(P, Q)$ is semi-harmonious. Then a simple use of (5.20) gives

$$\|PQ - P_1\| = \|(I-Q)(I-P) - P_4\| = \|[(I-Q)(I-P) - P_4]^*\|.$$

The equation above together with (6.5) yields the desired conclusion. $\square$

**Lemma 6.3.** *Suppose that $(P, Q)$ is a quasi-projection pair on $H$ such that $H_1$ defined by (3.4) is orthogonally complemented in $H$. Let $X$ be a Hilbert space and $\pi : \mathcal{L}(H) \to \mathbb{B}(X)$ be a unital $C^*$-morphism. Then*

$$\|\pi(P)\pi(Q) - \pi(P_{H_1})\| = \max\{\alpha, \|P_M - \pi(P_{H_1})\|\}, \qquad (6.6)$$
$$\|\pi(P)\pi(P_{\mathcal{R}(Q)}) - \pi(P_{H_1})\| = \max\{\beta, \|P_M - \pi(P_{H_1})\|\}, \qquad (6.7)$$

*where $M$, $\alpha$ and $\beta$ are given respectively by*

$$M = \mathcal{R}(\pi(P)) \cap \mathcal{R}(\pi(Q)), \quad \alpha = \|\pi(P)\pi(Q) - P_M\|, \qquad (6.8)$$
$$\beta = \|\pi(P)\pi(P_{\mathcal{R}(Q)}) - P_M\|.$$

*Proof.* Let $P_{H_1}$ be denoted simply by $P_1$, and let

$$S_1 = \pi(P)\pi(Q) - P_M, \quad S_2 = P_M - \pi(P_1). \qquad (6.9)$$

It is clear that $PP_1 = QP_1 = P_1$, so $\pi(P)\pi(P_1) = \pi(Q)\pi(P_1) = \pi(P_1)$, which implies that $\pi(P_1) \le P_M$, hence $S_2$ is a projection and $S_1 S_2 = 0$. It follows that
$$S_2 S_1^* = S_2^* S_1^* = (S_1 S_2)^* = 0.$$

Since $(\pi(P), \pi(Q))$ is also a quasi-projection pair, by (3.15) we have
$$\begin{aligned} M &= \mathcal{R}\big(\pi(P)\big) \cap \mathcal{R}\big(\pi(Q)^*\big) = \mathcal{R}\big(\pi(P)\big) \cap \mathcal{R}\big(\pi(Q^*)\big), \\ N &:= \mathcal{N}\big(\pi(P)\big) \cap \mathcal{N}\big(\pi(Q)\big) = \mathcal{N}\big(\pi(P)\big) \cap \mathcal{N}\big(\pi(Q)^*\big). \end{aligned} \qquad (6.10)$$

As $\pi(Q)^*$ is an idempotent and $M \subseteq \mathcal{R}\big(\pi(Q)^*\big)$, it gives $\pi(Q)^* P_M = P_M$. By (3.15) $H_1 \subseteq \mathcal{R}(Q^*)$, so
$$\pi(Q)^* \pi(P_1) = \pi(Q^* P_1) = \pi(P_1).$$

Therefore, it can be derived from (6.9) that $S_1^* S_2 = 0$, which yields $S_2 S_1 = (S_1^* S_2)^* = 0$. It follows that
$$\begin{aligned} &S_1 S_1^* \ge 0, \quad S_2 \ge 0, \quad S_1 S_2 = S_2 S_1^* = 0, \quad S_1 S_1^* \cdot S_2 = S_2 \cdot S_1 S_1^* = 0, \\ &(S_1 + S_2)(S_1 + S_2)^* = (S_1 + S_2)(S_1^* + S_2) = S_1 S_1^* + S_2. \end{aligned}$$

Since $S_2$ is a projection, we have $\|S_2\| = \|S_2\|^2$. Consequently,
$$\begin{aligned} \|S_1 + S_2\|^2 &= \|(S_1 + S_2)(S_1 + S_2)^*\| = \max\{\|S_1 S_1^*\|, \|S_2\|\} \\ &= \max\{\|S_1\|^2, \|S_2\|^2\}. \end{aligned}$$

This shows the validity of (6.6).

Following the same line as in the derivation of (6.6), we see that (6.7) is also true whenever $P_M \le \pi(P_{\mathcal{R}(Q)})$. So, it needs only to prove this latter inequality. It is known that $P_{\mathcal{R}(Q)} = Q(Q + Q^* - I)^{-1}$ (see e.g. [8, Theorem 1.3]). Since $\pi$ is a unital $C^*$-morphism, we have
$$\pi\big(P_{\mathcal{R}(Q)}\big) = \pi(Q)\big[\pi(Q) + \pi(Q)^* - I\big]^{-1} = P_{\mathcal{R}[\pi(Q)]}, \qquad (6.11)$$

which leads obviously to $P_M \le \pi\big(P_{\mathcal{R}(Q)}\big)$. $\square$

Now, we are in the position to provide the main result of this section.

**Theorem 6.4.** *Let $(P, Q)$ be a quasi-projection pair on $H$ such that $H_1$ and $H_4$ defined by (3.4) are orthogonally complemented in $H$. Then (6.3) is satisfied.*

*Proof.* We begin with the following simple observation:

$$b \le a \text{ and } \max\{b, c\} = \max\{b, d\} \implies \max\{a, c\} = \max\{a, d\}. \quad (6.12)$$

Indeed, utilizing $b \le a$ yields

$$\max\{a, c\} = \max\{a, b, c\} = \max\{a, \max\{b, c\}\} = \max\{a, \max\{b, d\}\}$$
$$= \max\{a, b, d\} = \max\{a, d\}.$$

Let $X$ be a Hilbert space and $\pi : \mathcal{L}(H) \to \mathbb{B}(X)$ be a faithful unital $C^*$-morphism. Denote $P_{H_1}$ and $P_{H_4}$ by $P_1$ and $P_4$, respectively. For simplicity, we use the same notation $I$ for the identity operators $I_H$ and $I_X$. Since $\big(\pi(P), \pi(Q)\big)$ is a quasi-projection pair on a Hilbert space, it is harmonious and hence is semi-harmonious. Therefore, a direct use of Lemma 6.2 yields

$$a := \|\pi(P)\pi(Q) - P_M\| = \|[I - \pi(P)][I - \pi(Q)] - P_N\|$$
$$= \|\pi(I - P)\pi(I - Q) - P_N\|, \quad (6.13)$$

where $M$ and $N$ are defined by (6.8) and (6.10), respectively. So, a simple use of Lemma 6.3 for quasi-projection pairs $(P, Q)$ and $(I - P, I - Q)$ gives

$$\|\pi(P)\pi(Q) - \pi(P_1)\| = \max\{a, c\}, \quad (6.14)$$
$$\|\pi(I - P)\pi(I - Q) - \pi(P_4)\| = \max\{a, d\}, \quad (6.15)$$

where $a$ is given by (6.13) and

$$c = \|P_M - \pi(P_1)\|, \quad d = \|P_N - \pi(P_4)\|.$$

Note that $\pi\big(P_{\mathcal{R}(Q)}\big)$ is a projection whose range equals $\mathcal{R}\big(\pi(Q)\big)$ (see (6.11)), so it is clear that

$$\big[\pi(P)\pi(Q) - P_M\big]\pi\big(P_{\mathcal{R}(Q)}\big) = \pi(P)\pi\big(P_{\mathcal{R}(Q)}\big) - P_M.$$

Hence
$$b := \big\|\pi(P)\pi\big(P_{\mathcal{R}(Q)}\big) - P_M\big\| \le a \cdot \big\|\pi\big(P_{\mathcal{R}(Q)}\big)\big\| \le a.$$

Now, we turn to prove that $\max\{b, c\} = \max\{b, d\}$. It is obvious that

$$\mathcal{R}(P) \cap \mathcal{R}(P_{\mathcal{R}(Q)}) = H_1, \quad \mathcal{N}\big(P_{\mathcal{R}(Q)}\big) = \mathcal{R}(Q)^\perp = \mathcal{N}(Q^*).$$

So, from (3.15) we can obtain

$$\mathcal{N}(P) \cap \mathcal{N}\big(P_{\mathcal{R}(Q)}\big) = H_4.$$

Moreover, due to (6.11) and (6.10) we have

$$\mathcal{R}(\pi(P)) \cap \mathcal{R}(\pi(P_{\mathcal{R}(Q)})) = \mathcal{R}(\pi(P)) \cap \mathcal{R}(\pi(Q)) = M,$$
$$\mathcal{N}(\pi(P)) \cap \mathcal{N}(\pi(P_{\mathcal{R}(Q)})) = \mathcal{N}(\pi(P)) \cap \mathcal{N}(\pi(Q)^*) = N.$$

Replacing $(P,Q)$ in Lemma 6.1 with $(P, P_{\mathcal{R}(Q)})$ and $(\pi(P), \pi(P_{\mathcal{R}(Q)}))$ respectively yields

$$\|PP_{\mathcal{R}(Q)} - P_1\| = \|(I-P)(I - P_{\mathcal{R}(Q)}) - P_4\|,$$
$$b = \|(I - \pi(P))(I - \pi(P_{\mathcal{R}(Q)})) - P_N\|.$$

We may combine Lemma 6.3 with the above two equations to conclude that

$$\max\{b,c\} = \|\pi(PP_{\mathcal{R}(Q)} - P_1)\| = \|\pi[(I-P)(I - P_{\mathcal{R}(Q)}) - P_4]\|$$
$$= \|[I - \pi(P)][I - \pi(P_{\mathcal{R}(Q)})] - \pi(P_4)\| = \max\{b,d\}.$$

It follows from (6.12) that $\max\{a,c\} = \max\{a,d\}$. Hence, (6.3) can be derived from (6.14) and (6.15) together with the faithfulness of $\pi$. □